\documentclass[final]{siamltex}

\usepackage{amsfonts}
\usepackage{amssymb}
\usepackage{amsbsy}

\usepackage{pst-grad, pst-node}
\usepackage{pstricks,pst-plot,pst-eps}
\usepackage{epsfig}
\usepackage{verbatim}
\usepackage{color}
\usepackage{cases}

\def\Section{\setcounter{equation}{0} \setcounter{thm}{0}\section}
\def\thesection {\arabic{section}}
\def\thesubsection {\thesection.\arabic{subsection}}
\def\thesubsubsection {\thesubsection.\arabic{subsubsection}}
\def\theequation {\thesection.\arabic{equation}}
\def\be{\begin{equation}\displaystyle}
\def\ee{\end{equation}}
\def\bel{\begin{equation} \displaystyle \begin{array}{l} }
\def\eel{\end{array} \end{equation} }
\def\bell{\begin{equation} \displaystyle \begin{array}{ll}  }
\def\eell{\end{array} \end{equation} }

\def\bea{\begin{eqnarray}}
\def\eea{\end{eqnarray} }
\def\bean{\begin{eqnarray*}}
\def\eean{\end{eqnarray*} }

\def\CC{\mathbb{C}}

\def\RR{\mathbb{R}}

\def\O{\mathcal{O}}

\catcode`@=11
\def\eqalign#1{\null\,\vcenter{\openup1\jot \m@th
   \ialign{\strut \hfil$\displaystyle{##}$ & $\displaystyle{{}##}$\hfil
      \crcr#1\crcr}}\,}
\catcode`@=12

\newtheorem{thm}{Theorem}  
\newtheorem{lem}[thm]{Lemma}
\newtheorem{cor}[thm]{Corollary}

\newtheorem{defi}[thm]{Definition}
\newtheorem{prop}[thm]{Proposition}
\newtheorem{rem}[thm]{Remark}

\def\debproof {\noindent{\em Proof}. }
\def\finproof{\hfill \endproof}  

\def\debthm {\begin{thm}}
\def\finthm {\end{thm}}
\def\deblem {\begin{lem}}
\def\finlem {\end{lem}}
\def\debprop {\begin{prop}}
\def\finprop {\end{prop}}
\def\debcor {\begin{cor}}
\def\fincor {\end{cor}}
\def\debdef {\begin{defi}}
\def\findef {\end{defi}}
\def\debrem {\begin{rem}}
\def\finrem {\end{rem}}

\def\debbox{ 
   \hspace{10mm}\smallskip\newline 
   \hspace*{3mm}\begin{tabular}[b]{|r |l} \hline \\
   \rm\makebox[2mm]{}\begin{minipage}{\noulong}}
\def\finbox{\end{minipage}\\ \\ \hline \end{tabular}}

\psset{xunit=0.4cm,yunit=0.4cm} 
\definecolor{dgreen}{rgb}{0,.8,0.2}

\def\ds{\displaystyle}

\def\eps{\varepsilon}

\def\bar#1{{\overline #1}}

\def\i{{i}} 

\newcommand{\eqref}[1]{(\ref{#1})}

\catcode`@=11
\def\supess{\mathop{\operator@font ess\,sup}}
\catcode`@=12

\title{Stationary Schr\"odinger equation in
      the semi-classical limit: numerical coupling of oscillatory and evanescent regions}

\author{Anton Arnold\thanks{Institut f\"ur Analysis und Scientific Computing, Technische Universit\"at Wien, Wiedner Hauptstr. 8, A-1040 Wien, Austria ({\tt anton.arnold@tuwien.ac.at}).}
\and Claudia Negulescu\thanks{Institut de Math\'ematiques de Toulouse (UMR 5219), Universit\'e de Toulouse, CNRS; UPS IMT, F-31062 Toulouse Cedex 9, France.}  ({\tt claudia.negulescu@math.univ-toulouse.fr}).}

\begin{document}
\maketitle

\begin{abstract}
This paper is concerned with a 1D Schr\"odinger scattering problem involving both oscillatory and evanescent regimes, separated by jump discontinuities in the potential function, to avoid ``turning points''. We derive a non-overlapping domain decomposition method to split the original problem into sub-problems on these regions, both for the continuous and afterwards for the discrete problem. Further, a hybrid WKB-based numerical method is designed for its efficient and accurate solution in the semi-classical limit: a WKB-marching method for the oscillatory regions and a FEM with WKB-basis functions in the evanescent regions. We provide a complete error analysis of this hybrid method and illustrate our convergence results by numerical tests.


\end{abstract}



\pagestyle{myheadings}
\thispagestyle{plain}
\markboth{A. ARNOLD AND C. NEGULESCU}{DOMAIN COUPLING FOR THE 1D SCHR\"ODINGER EQUATION
}


\Section{Introduction}\label{SEC1}
This paper deals with the design, error analysis, and numerical study of an asymptotically correct scheme for the  numerical solution of the stationary Schr\"odinger equation in one dimensional scattering situations:
\be \label{EQ_ref}
\eps^2 \psi''(x) + a(x) \psi(x) =0\,, \quad x \in \RR\,,
\ee
where $0<\eps \ll 1$ is a very small parameter and $a(x)$ a piecewise (sufficiently)
smooth, real function. On the one hand, for $a(x)>0$, the solution is highly oscillatory, with the small (local) wave length $\lambda(x)={2\pi\eps\over \sqrt{a(x)}}$. On the other hand, for $a(x)<0$, the wave function $\psi$ is (essentially) exponentially de/increasing, typically exhibiting a thin boundary layer with thickness of the order  $\O\big({\eps\over \sqrt{|a(x)|}}\big)$. A key aspect of this paper is that $a(x)$ takes both signs. Hence, we have to cope with a classical multi-scale problem, combining different types of arduousnesses and multi-scale behaviours. Numerically, we aim at recovering these fine structures of the solution, however \emph{without} using a fine spatial grid. To this end we shall develop a (non-overlapping) domain decomposition method (DDM) to separate the oscillatory and evanescent regions, as they require very different numerical approaches. This DDM allows to recover at the continuous level the exact analytical solution in a single sweep (ag!
 ainst the direction of the 
incoming plane wave) with appropriate interface conditions and a final scaling.\\

The study of such multi-scale problems is very challenging from a theoretical as well as numerical point of view. In both situations (or regions) a classical discretization (like in \cite{Ba1,Ba2}) 
requires a very  fine mesh in order to accurately resolve 
the oscillations 
and boundary layers.
However, with a step size requirement of $h={\cal O}(\eps)$,  
standard numerical methods would be very costly and inefficient here. 

\begin{figure}[htbp]
\begin{center}
\includegraphics[width=6cm]{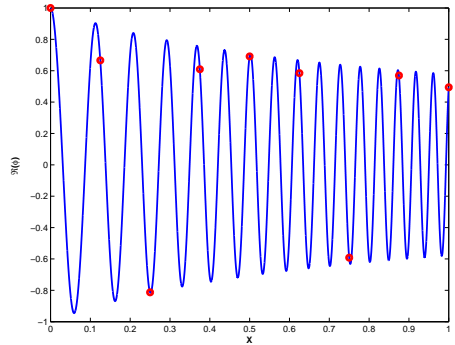}
\end{center}
\caption{\label{onde} {\footnotesize In standard numerical methods
highly oscillating solutions require a very fine mesh to capture the oscillations.
However, with the analytic pre-processing of our WKB-marching method an accurate solution can be obtained on a coarse grid
(dots). Plotted is the solution $\Re\psi(x)$ of \eqref{EQ_ref} with $\eps=0.01,\,h=0.125,$ and $a=(x+\frac12)^2$.}}
\end{figure}

Concerning the \emph{oscillatory case}, several WKB-based numerical schemes (named after the physicists Wentzel, Kramers, and Brillouin) have been developed and analysed in the last decade. Their  goal is to use a coarse spatial grid with step size $h>\lambda$ (see Figure \ref{onde}), reducing the limitation to at least $h={\cal O}(\sqrt\eps)$. For marching schemes we refer to \cite{ABAN11, JL03, lub}, whereas a finite element method (FEM) using oscillatory WKB-basis functions was introduced in \S2 of \cite{Ne05} and in \cite{Cla}. This FEM approach has the disadvantage that it requires a non-resonance condition between the mesh-size $h$ and the wave length $\lambda$ of the solution. By contrast, this restriction is not necessary in the above mentioned marching schemes.

Numerical approaches for the \emph{evanescent regime} (as $\eps\to0$) have been considered much less, so far. We refer to \S3 in \cite{Ne05} for the formulation of a WKB-based multiscale FEM-scheme, including its numerical coupling to the oscillatory region (also based on a FEM). But a numerical analysis has, to our knowledge, not been carried out yet. This paper also aims at closing this gap.
In this evanescent regime the problem \eqref{EQ_ref} is elliptic, and for the example of $a=const.$, a solution  is given by a linear combination of the basis-functions $\exp(\pm\frac{\sqrt{|a|}}{\eps}x)$. Hence this region must be considered as a boundary value problem (BVP) and solved e.g.\ by a finite difference or a FEM method. A reformulation as an initial value problem (IVP) and the use of a marching scheme would be inherently unstable, due to the unbounded growth of $\exp(\frac{\sqrt{|a|}}{\eps}x)$ (in $\eps$). This contrasts with the oscillatory regime, where a basis of the solution  would be given by the bounded functions $\exp(\pm i\frac{\sqrt{a}}{\eps}x)$.
Consequently, we are faced with coupling two different approaches for solving \eqref{EQ_ref} in the semi-classical limit: an IVP for $a(x)>0$ and a BVP for $a(x)<0$.\\

The goal of this paper is to analyse the numerical coupling of oscillatory and evanescent regimes, using WKB-ansatz functions for both situations. In the oscillatory regime we shall use the marching scheme from \cite{ABAN11}, and in the evanescent regime a FEM like in \S3 of \cite{Ne05}. In the first case, the key idea is to eliminate analytically the dominant oscillations of the solution to \eqref{EQ_ref}. The transformed problem then has a much smoother solution, in the sense that the amplitude of the residual oscillations is much smaller than in the original problem -- often by many orders of magnitude (in fact by the order $\eps^2$, cf. Propositions 2.1, 2.2 in \cite{ABAN11}). Hence, the new problem can be solved numerically on a coarse grid, still yielding a very accurate approximation.
In the evanescent regime, the key idea of the WKB-FEM method is to use WKB-ansatz functions (of exponential type), rather than the standard polynomials. Finally the strategy is the same as in the oscillatory region: to filter out the boundary layer via well-chosen basis functions. Since WKB-basis functions are asymptotic solutions to \eqref{EQ_ref} (as $\eps\to0$), this method is again very accurate on a coarse grid. In this paper we shall provide first the numerical analysis of the WKB-FEM method for the evanescent regime (from scratch), and then the error analysis of the hybrid DDM -- building on the convergence results in \cite{ABAN11} for the WKB-marching method.\\

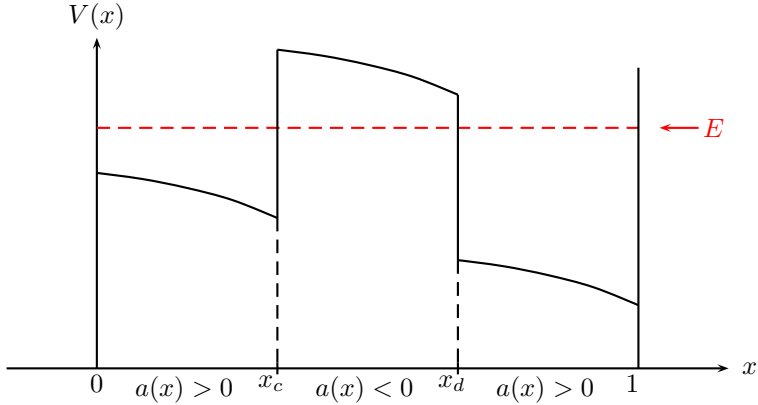
\begin{figure}[tb]
\begin{center}
\begin{pspicture}
(-4,-2)(23,15)
\psline{->}(-3,0)(21,0)  
\psline{->}(0,0)(0,11)   
\rput(21.7,0){$x$}
\rput(0,11.7){$V(x)$}
\psline[linecolor=red,linestyle=dashed]{-}(0,8)(18,8)
\psline[linecolor=red]{<-}(18.7,8)(20,8)
\rput(20.5,8){$\color{red} E$}
\pscurve(0,6.5)(1.5,6.3)(3,6.0)(4.5,5.6)(6,5)
\psline(6,5)(6,10.6)
\pscurve(6,10.6)(7.5,10.4)(9,10.1)(10.5,9.7)(12,9.1)
\psline(12,9.1)(12,3.6)
\pscurve(12,3.6)(13.5,3.4)(15,3.1)(16.5,2.7)(18,2.1)
\psline(18,0)(18,10)
\psline(6,-0.2)(6,0.2)
\psline[linestyle=dashed](6,0.2)(6,5)
\psline(12,-0.2)(12,0.2)
\psline[linestyle=dashed](12,0.2)(12,3.6)
\rput(2.9,-0.7){$a(x)>0$}
\rput(8.9,-0.7){$a(x)<0$}
\rput(14.9,-0.7){$a(x)>0$}
\rput(0,-0.5){0}
\rput(5.8,-0.5){$x_c$}
\rput(11.8,-0.5){$x_d$}
\rput(17.8,-0.5){$1$}

\end{pspicture}
\end{center}
\caption{\label{injection} {\footnotesize Sketch of a tunnelling structure with two oscillatory and one evanescent regions. Electrons are injected from the right boundary with energy $E$. The coefficient function in \eqref{EQ_ref} is $a(x):=E-V(x)$.}}
\end{figure}

Problems similar to  \eqref{EQ_ref}  or in general that require the numerical integration of highly oscillatory equations play an essential role in a wide range of
physical phenomena: e.g. electromagnetic and acoustic scattering
(Maxwell and Helmholtz equations in the high frequency regime), wave evolution problems in quantum and
plasma physics (Schr\"odinger equation in the semi-classical regime), and stiff mechanical systems. The application we are interested in here, stems from the electron transport in 
nano\-scale semiconductor devices, like quantum wave-guides \cite{Ar08},
resonant tunnelling diodes (RTDs) \cite{BAP06}, MOSFETs \cite{Cla_Naoufel}, etc. 
In a 1D model setting, which is appropriate for RTDs or for the longitudinal dynamics in each transversal mode in MOSFETs,
the governing equation is the stationary Schr\"odinger equation.
In an idealized model we assume that the given electrostatic potential $V(x)$ is constant in the left lead $x\in(-\infty,0]$, with value $V(0)$, and equally in the right lead $x\in[1,\infty)$, with value $V(1)$. Hence the Schr\"odinger equation can be 
complemented with open boundary conditions at both ends:
\be \label{SchBVP}
\left\{
\begin{array}{l}
\ds - \eps^2 \psi ''(x) + V(x) \psi(x)=E \psi(x)\,, \quad x \in (0,1)\,, \\[3mm]
\ds \eps\psi'(0)+{i} \sqrt{a(0)} \,\psi(0) =0\,, \quad\mbox{ if }a(0)>0\,,\\[3mm]
\ds \eps\psi'(0)- \sqrt{|a(0)|} \,\psi(0) =0\,, \quad\mbox{ if }a(0)<0\,,\\[3mm]
\ds \eps\psi'(1) -{i} \sqrt{a(1)} \,\psi(1) = -2 {i} \sqrt{a(1)}\,.
\end{array}
\right.
\quad a(x):=E-V(x)\,,
\ee
This equation describes the state of an electron that is injected with the prescribed energy $E>0$ from the right boundary (or lead)
into an electronic device (diode, e.g.), modelled on the interval $[0,1]$. The corresponding (complex valued) wave function is denoted by $\psi(x)$, 
where $|\psi(x)|^2$ is related to the spatial probability density of the electron. Due to the continuous injection of a plane wave function at $x=1$, we cannot expect $|\psi|^2$ to be normalised here (in $L^1(0,1)$).
When considering the equivalent problem on $\RR$, $\psi$ rather describes a scattering state with $\psi\in L^\infty(\RR)$.
The small parameter $0<\eps<1$ is the re-scaled Planck constant.
To make the link with \eqref{EQ_ref}, the coefficient function $a(x)$ is given by $a(x):=E-V(x)$.
To allow for an injection at $x=1$, we have to require that $a(1)>0$, cf.\ Figure \ref{injection}. In fact, $E>V(x)$ characterises the oscillatory, classically allowed regime, whereas $E<V(x)$ characterises the evanescent, classically forbidden regime. Fig.\ \ref{injection} sketches a tunnelling structure including both regimes, which are rather different. The boundary conditions in \eqref{SchBVP} are the so called open or \emph{transparent boundary conditions}, permitting an electron wave to
enter or leave the device without reflections \cite{lent}. Due to the injected plane wave of electrons, the boundary condition at $x=1$ is inhomogeneous. But at $x=0$ it is homogeneous, due to the free outflow of the electron wave.

In the present work we shall not discuss (in detail) situations with \emph{turning points}, i.e.\ zeros of $a$, but rather concentrate on devices with abrupt jumps at the interfaces between oscillatory and evanescent regions. This first step is simpler to treat and will be extended in a subsequent work. In Section \ref{SEC5} we shall comment on situations incorporating turning points.

\smallskip
For the solvability of this model, the following simple result holds:
\vspace{0.2cm}
\begin{prop}\label{Prop1.1}
Let $a\in L^\infty(0,1)$ with $a(0)\ne0$ and $a(1)>0$.\footnote{Here and in the sequel, $a(0)$ and $a(1)$ are not meant as the point values of the function $a$ (which would not be defined for an $L^\infty$--function), but rather as the constant potential in the left and, resp., right leads.}
Then the boundary value problem \eqref{SchBVP} has a unique solution $\psi\in W^{2,\infty}(0,1)\hookrightarrow C^1[0,1]$.
\end{prop}

\vspace{0.2cm}
\debproof
For the case of an oscillatory outflow, i.e. $a(0)>0$, the proof was provided in Proposition 2.3 of \cite{NDM}. For an evanescent outflow, i.e. $a(0)<0$, the proof is analogous (multiplying the Schr\"odinger equation by $\bar \psi$, integrating by parts, and taking the imaginary part).
\finproof
\medskip

\noindent {\bf WKB-technique.}\\
Both parts of the hybrid numerical method studied in \S \ref{SEC3} will be based on WKB functions. Hence, let us first review the well-known WKB-approximation (cf.\ \cite{LL85}) for the singularly perturbed ODE \eqref{EQ_ref}. In the standard approach, for the oscillatory regime (i.e.\ $a>0$), the WKB-ansatz
\begin{equation}\label{WKB}
  \psi(x)\sim\exp\left(\frac1\eps \sum_{p=0}^\infty \eps^p \phi_p(x)\right)\,,
\end{equation}
is inserted in \eqref{EQ_ref} and after comparison of the
$\eps^p$-terms, leads to
\begin{eqnarray}\label{WKB-fct}
\phi_0(x)&=&\pm \i\int_0^x\sqrt{a(y)}\,dy {+const.}\,,\nonumber\\
\phi_1(x)&=&\ln a(x)^{-1/4}+const.\,,\\
\phi_2(x)&=&\mp \i \int_0^x\beta(y)\,dy {+const.}\,,\quad  \beta:=
-\frac{1}{2|a|^{1/4}}(|a|^{-1/4})''\,.
\nonumber
\end{eqnarray}
Truncating the ansatz \eqref{WKB} after $p=0,\,1,$ or 2, yields the asymptotic approximation 
$\psi(x)\approx C\varphi_p^{os}(x)$, with the following \emph{oscillatory WKB-functions} (of the three lowest orders in $\eps$):
\begin{eqnarray}\label{osWKB-fct}
  \varphi_0^{os}(x)&=&\exp\left(\pm\frac{\i}{\eps}\int_0^x\sqrt{a(y)}\,dy  \right)\,,\nonumber\\
  \varphi_1^{os}(x)&=&\frac{\exp\left(\pm\frac{\i}{\eps}\int_0^x\sqrt{a(y)}\,dy  \right)}{\sqrt[4]{a(x)}}\,,\\
  \varphi_2^{os}(x)&=&\frac{\exp\left(\pm\frac{\i}{\eps}\int_0^x\big[\sqrt{a(y)}-\eps^2\beta(y)\big]\,dy  \right)}{\sqrt[4]{a(x)}}\,.\nonumber
\end{eqnarray}

Proceeding analogously for the evanescent regime (i.e.\ $a<0$) yields the following \emph{evanescent WKB-functions} (of the three lowest orders in $\eps$):
\begin{eqnarray}\label{evWKB-fct}
  \varphi_0^{ev}(x)&=&\exp\left(\pm\frac{1}{\eps}\int_0^x\sqrt{|a(y)|}\,dy  \right)\,,\nonumber\\
  \varphi_1^{ev}(x)&=&\frac{\exp\left(\pm\frac{1}{\eps}\int_0^x\sqrt{|a(y)|}\,dy  \right)}{\sqrt[4]{|a(x)|}}\,,\\
  \varphi_2^{ev}(x)&=&\frac{\exp\left(\pm\frac{1}{\eps}\int_0^x\big[\sqrt{|a(y)|}+\eps^2\beta(y)\big]\,dy  \right)}{\sqrt[4]{|a(x)|}}\,.\nonumber
\end{eqnarray}

In the hybrid numerical method analysed in \S \ref{SEC3} we shall use the first order WKB ansatz functions $\varphi_1^{ev}$ for the FEM in the evanescent region. And in the oscillatory region we shall use the second order WKB functions $\varphi_2^{os}$ to transform \eqref{EQ_ref} into a smoother problem that can be solved accurately and efficiently on a coarse grid. Since we shall use two different numerical approaches in the two regimes, also the corresponding error orders will be rather different (both with respect to $\eps$ and to the grid size $h$). Hence there is,   a-priori, no obvious natural choice for the orders of the two WKB-methods. We choose a first order WKB-method in the evanescent region to keep 
the complexity of the numerical scheme and 
the technicalities of its convergence analysis low. Furthermore we choose a second order WKB-method in the oscillatory region such that we can use the results from \cite{ABAN11} (without having to redo a first order WKB-analysis).
Anyhow, our hybrid method is second order with respect to $h$.\\

This paper is organised as follows. In Section \ref{SEC2} we present and analyse the (non-overlapping) domain decomposition method for the singularly perturbed ODE (\ref{EQ_ref}) on the continuous level. Propositions \ref{prop2.2} and \ref{prop2.5} establish that this DDM yields the analytical solution in one sweep for cases consisting of two and, resp., three distinct regions. In \S\ref{SEC3} we first review the two different numerical WKB-methods for the two distinct regions and establish convergence of the WKB-FEM. Then we prove convergence of the overall hybrid WKB-method (WKB-FEM in the evanescent regime coupled to a WKB-marching scheme for the oscillatory region), with Theorem \ref{th:hybrid-conv} as the main result. 
In \S\ref{SEC4} we illustrate our convergence results on some numerical examples treated with our scheme, including an example with a tunnelling structure. 
And finally, in \S\ref{SEC5} we briefly discuss extensions of our WKB method to coefficient functions with turning points.

\Section{Domain decomposition of the Schr\"odinger boundary value problem} \label{SEC2}
In this section we shall consider the Schr\"odinger BVP \eqref{SchBVP} with given coefficient functions $a(x)=E-V(x)$ corresponding to two different scenarios -- first two coupled regions, then three regions. We confine ourselves here to these examples only for practical reasons: This setup already shows all the interesting features of the BVP, and it can be generalised easily to more regions.

\subsection{Two coupled regions}\label{SEC2.1}

\begin{figure}[tb]
\begin{center}
\begin{pspicture}
(-4,-2)(23,15)
\psline{->}(-3,0)(21,0)  
\psline{->}(0,0)(0,12)   
\rput(21.7,0){$x$}
\rput(0,12.7){$V(x)$}
\psline[linecolor=red,linestyle=dashed]{-}(0,8)(18,8)
\psline[linecolor=red]{<-}(18.7,8)(20,8)
\rput(20.5,8){$\color{red} E$}
\pscurve(0,10)(2,10.6)(4,10.4)(6,9.7)(8,9.5)(10,9.1)
\psline(10,9.1)(10,3.6)
\pscurve(10,3.6)(12,3.5)(13.5,3.4)(15,3.1)(16.5,2.7)(18,2.1)
\psline(18,0)(18,10)
\psline[linecolor=dgreen]{->}(-1,9)(1.5,9)
\psline[linecolor=dgreen]{->}(12.75,9)(15.25,9)
\rput(1.5,9.8){$\color{dgreen} solve$}
\psline(10,-0.2)(10,0.2)
\psline[linestyle=dashed](10,0.2)(10,3.6)
\rput(4.9,-0.7){$a(x)<0$}
\rput(13.9,-0.7){$a(x)>0$}
\rput(0,-0.5){0}
\rput(9.9,-0.5){$x_d$}
\rput(17.9,-0.5){$1$}
\pscircle(5,0.7){0.25}
\rput(5.0,0.7){$1$}
\pscircle(14,0.7){0.25}
\rput(14.0,0.7){$2$}
\rput(5,7){$\psi(x)=\alpha\chi(x)$}
\rput(14,7){$\psi(x)=\alpha\varphi(x)$}

\end{pspicture}
\end{center}
\caption{\label{2regions} {\footnotesize Potential barrier: While electrons are injected from the right boundary with energy $E>0$,
the decomposed problem has to be solved from left to right (as a BVP--IVP). }}
\end{figure}
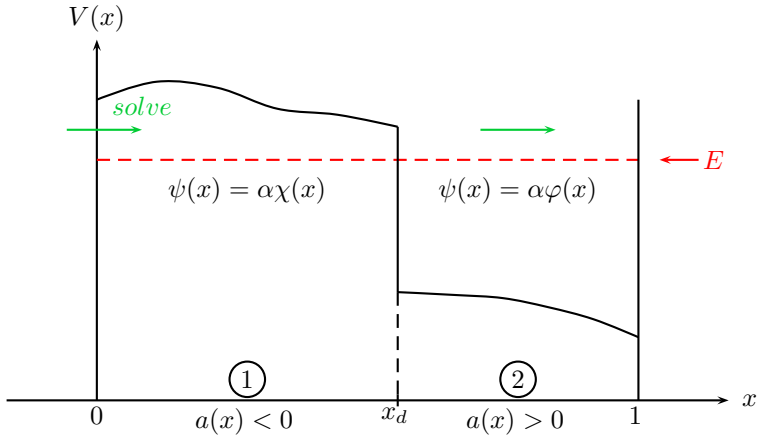

We start with the situation illustrated in Figure \ref{2regions}: It consists of two regimes, an evanescent region with $a:=E-V<0$ (adjacent to the left boundary) and an oscillatory region with $a>0$ (adjacent to the right boundary). Since we exclude turning points here, the function $a$ is assumed to have a jump discontinuity (and a sign change) at the interface $x=x_d$. 
Moreover, for this section we shall assume:\\

{\bf Hypothesis a2} {\it Let $a\in L^\infty(0,1)$ with $a\big|_{(0,x_d)}<0$,  $a\big|_{(x_d,1)}>0$, $a(0)<0$, and $a(1)>0$.}\\

Following the basic idea from \cite{ABAN11} we shall solve the BVP \eqref{SchBVP} as two consecutive sub-problems: We start with the evanescent region $(0,x_d)$, where a BVP is solved (for stability reasons, as mentioned in \S\ref{SEC1}). This is followed by an IVP on the oscillatory region $(x_d,1)$. So we shall proceed in the opposite direction of the injection direction (see Figure \ref{2regions}). Due to Proposition \ref{Prop1.1}, the solution to \eqref{SchBVP} is in $C^1[0,1]$. Hence the solutions on these two sub-regions are matched by continuity of $\psi$ and $\psi'$ at $x=x_d$.\\

For the BVP on $(0,x_d)$, the original problem \eqref{SchBVP} provides only one homogeneous Robin boundary condition (BC) at $x=0$. Hence, we supply the BVP with an auxiliary, artificial BC at $x=x_d$. Here, both an inhomogeneous Dirichlet or Neumann BC would work from an algorithmic point of view. In order to simplify the numerical analysis in \S\ref{SEC311} below, we shall use at this point $\eps\chi'(x_d)=1$ for the auxiliary wave function $\chi$. 
While this auxiliary value has the correct $\eps$--order (cf. Proposition \ref{prop2.2a} and Lemma \ref{U-estimate} below), it
will in general not be the correct derivative of the global solution $\psi$. 
Its correct value will finally be obtained by scaling the auxiliary functions using the remaining inhomogeneous Robin BC at $x=1$ (cf. \eqref{SchBVP}). This leads to the following domain decomposition and problem coupling for the auxiliary wave functions $\chi,\,\varphi$:\\

\noindent
\underline{Step 1 -- BVP for $\chi$ in region (1):}
\be \label{region1_ev}
\,\,\, \left\{
\begin{array}{l}
\ds \eps^2 \chi ''(x) + a(x) \chi(x)=0\,, \qquad x \in (0,x_d)\,,  \qquad \quad a(x)= E -V(x) <0\,,\\[3mm]
\ds \eps\, \chi'(0) - \sqrt{|a(0)|}\, \chi(0)=0\,, \quad \mbox{(Robin BC for $\chi$ )}\\[3mm]
\ds \eps \chi'(x_d) =1\,. \qquad\qquad\qquad\quad\: \mbox{(auxiliary Neumann BC )}
\end{array}
\right.
\ee
\medskip

\noindent
\underline{Step 2 -- IVP for $\varphi$ in region (2):}
\be \label{region2_ev}
\,\,\, \left\{
\begin{array}{l}
\ds \eps^2 \varphi ''(x) + a(x) \varphi(x)=0\,, \qquad x \in (x_d,1)\,, \qquad \quad a(x)= E -V(x)>0\,, \\[3mm]
\ds \varphi(x_d) =\chi(x_d)\,, \qquad\qquad\:\: \mbox{(implies continuity of $\psi$ at $x_d$)}\\[3mm]
\ds \eps\varphi'(x_d) =\eps\chi'(x_d)=1 \,.\quad\: \mbox{(implies continuity of $\psi'$ at $x_d$)}
\end{array}
\right.
\ee
\medskip

\noindent
\underline{Step 3 -- Scaling of the auxiliary wave functions:}
\be \label{psi-scaling_ev}
\psi(x):=\left\{
\begin{array}{l}
\ds \alpha\,\chi(x)\,, \quad x \in (0,x_d)\,, \\[3mm]
\ds \alpha\,\varphi(x)\,, \quad x \in (x_d,1)\,,
\end{array}
\right.
\ee
with the scaling parameter $\alpha\in\CC$ defined via
\be\label{alpha-eq_ev}
   \qquad\alpha=\alpha(\varphi(1),\,\varphi'(1))=\frac{-2\i \sqrt{a(1)}}{\eps\varphi'(1)-i \sqrt{a(1)}\,\varphi(1)} \,, \quad \mbox{(due to the right BC in \eqref{SchBVP})}\,.
\ee
We note that the denominator in this expression for $\alpha$ is non-zero: On the one hand $\chi$ and $\varphi$ are both real valued, and on the other hand $\varphi(1)$ and $\varphi'(1)$ cannot vanish simultaneously (as otherwise $\varphi\equiv0$ would contradict the Neumann BC in \eqref{region1_ev}).

In this whole section we shall only require that $a\in L^\infty(0,1)$. As in Proposition \ref{Prop1.1}, $a(0)$ and $a(1)$ are hence not meant as the point values of the function $a$, 
but rather as the constant potential in the left and right leads.\\

Next we address the solvability of the algorithm  \eqref{region1_ev}-\eqref{alpha-eq_ev}.
\vspace{0.2cm}
\begin{lem} \label{PB-2zones}
Let Hypothesis a2 be satisfied.
Then:
\begin{enumerate}
\item[a)] The BVP \eqref{region1_ev} has $\forall \eps >0$ a unique solution $\chi\in W^{2,\infty}(0,x_d) \subset C^1[0,x_d]$.
\item[b)] The IVP \eqref{region2_ev} has $\forall \eps >0$ a unique solution $\varphi\in W^{2,\infty}(x_d,1)\subset C^1[x_d,1]$.
\end{enumerate}
Both functions $\chi$ and $\varphi$ are real functions and the parameter $\alpha\in \CC$ given by \eqref{alpha-eq_ev} is well-defined.
\end{lem}\\
\vspace{0.2cm}

The above lemma, whose proof is very easy, shows that the \emph{domain decomposition algorithm} \eqref{region1_ev}-\eqref{alpha-eq_ev} yields a unique function $\psi$ that is piecewise in $W^{2,\infty}$ and piecewise (on the two regions) a solution to the Schr\"odinger equation \eqref{SchBVP}. In fact, this DDM yields the unique solution of \eqref{SchBVP} as stated in the following proposition:
\vspace{0.2cm}
\debprop\label{prop2.2}
Let Hypothesis a2 be satisfied.
Then the function $\psi$ obtained from \eqref{region1_ev}-\eqref{alpha-eq_ev} belongs to $W^{2,\infty}(0,1)$ and is the unique solution of \eqref{SchBVP} (as guaranteed by Proposition \ref{Prop1.1}).
\finprop
\vspace{0.2cm}

\debproof
  The initial conditions in \eqref{region2_ev} imply $C^1$--continuity of $\psi$ at $x_d$. Hence, $\psi\in C^1[0,1]$, and this proves the claim.
\finproof

The following result provides the uniform-in-$\eps$ boundedness of this solution $\psi$. It generalizes Theorem 2.2 from \cite{Cla}, which holds only for one purely oscillatory region:
\vspace{0.2cm}
\debprop\label{prop2.2a}
Let Hypothesis a2 hold. Moreover, let the potential in the oscillatory region satisfy $a\in W^{1,\infty}(x_d,1)$ and $0<\tau_{os}\le a(x)$ $\forall x\in(x_d,1)$.
Then, the solution of \eqref{SchBVP} satisfies
\begin{equation}\label{2-zone-bound}
  \|\psi\|_{L^\infty(0,1)} + \eps \|\psi'\|_{L^\infty(0,1)} \le C\,,
\end{equation}
independently of $0<\eps\le1$.
\finprop

\vspace{0.2cm}
The simple, but lengthy proof is deferred to the Appendix.


\subsection{Three coupled regions}\label{SEC22}

\begin{figure}[tb]
\begin{center}
\begin{pspicture}
(-4,-2)(23,15)
\psline{->}(-3,0)(21,0)  
\psline{->}(0,0)(0,11)   
\rput(21.7,0){$x$}
\rput(0,11.7){$V(x)$}
\psline[linecolor=red,linestyle=dashed]{-}(0,8)(18,8)
\psline[linecolor=red]{<-}(18.7,8)(20,8)
\rput(20.5,8){$\color{red} E$}
\pscurve(0,6.5)(1.5,6.3)(3,6.0)(4.5,5.6)(6,5)
\psline(6,5)(6,10.6)
\pscurve(6,10.6)(7.5,10.4)(9,10.1)(10.5,9.7)(12,9.1)
\psline(12,9.1)(12,3.6)
\pscurve(12,3.6)(13.5,3.4)(15,3.1)(16.5,2.7)(18,2.1)
\psline(18,0)(18,10)
\psline[linecolor=dgreen]{->}(-1,9)(1.5,9)
\psline[linecolor=dgreen]{->}(13.75,9)(16.25,9)
\rput(1.5,9.8){$\color{dgreen} solve$}
\psline(6,-0.2)(6,0.2)
\psline[linestyle=dashed](6,0.2)(6,5)
\psline(12,-0.2)(12,0.2)
\psline[linestyle=dashed](12,0.2)(12,3.6)
\rput(2.9,-0.7){$a(x)>0$}
\rput(8.9,-0.7){$a(x)<0$}
\rput(14.9,-0.7){$a(x)>0$}
\rput(0,-0.5){0}
\rput(5.8,-0.5){$x_c$}
\rput(11.8,-0.5){$x_d$}
\rput(17.8,-0.5){$1$}
\pscircle(3,0.7){0.25}
\rput(3.0,0.7){$1$}
\pscircle(9,0.7){0.25}
\rput(9.0,0.7){$2$}
\pscircle(15,0.7){0.25}
\rput(15.0,0.7){$3$}
\rput(3,7){$\psi(x)=\beta\zeta(x)$}
\rput(9,7){$\psi(x)=\alpha\chi(x)$}
\rput(15,7){$\psi(x)=\alpha\varphi(x)$}

\end{pspicture}
\end{center}
\caption{\label{3regions} {\footnotesize Tunnelling structure: While electrons are injected from the right boundary with energy $E$,
the decomposed problem has to be solved from left to right (as an IVP--BVP--IVP). The coefficient function in \eqref{EQ_ref} is $a(x):=E-V(x)$.}}
\end{figure}
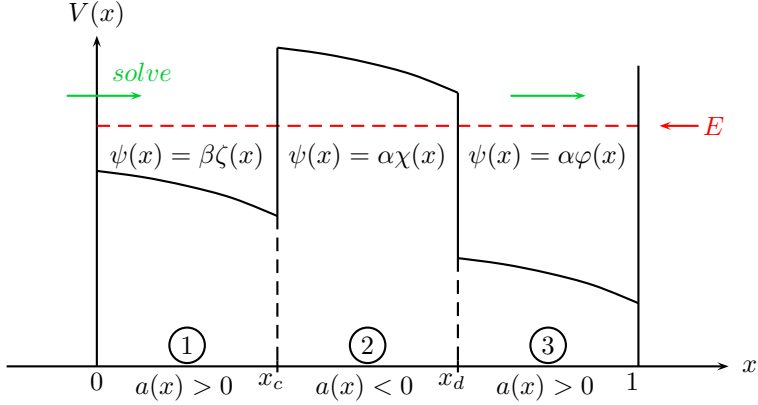

In this subsection we consider the Schr\"odinger BVP \eqref{SchBVP} with a given coefficient function $a(x)$ as illustrated in Figure \ref{3regions}: It consists of three regimes, two oscillatory regions at the interval boundaries and an evanescent region in the middle. Since we exclude turning points here, $a$ is assumed to have jump discontinuities (and sign changes) at the interfaces $x=x_c$ and $x=x_d$. The solution $\psi$ to the BVP \eqref{SchBVP} for such an example is illustrated in Fig.\  \ref{fig:4.1} below.
Moreover, for this section we shall assume on $a(x)$:\\

{\bf Hypothesis a3} {\it Let $a\in L^\infty(0,1)$ with $a\big|_{(x_c,x_d)}<0$,  $a\big|_{(0,x_c)\cup(x_d,1)}>0$, $a(0)>0$, and $a(1)>0$.}\\

Following the basic idea from \cite{ABAN11} we shall solve the BVP \eqref{SchBVP} as an IVP-BVP-IVP problem in the opposite direction of the injection direction, i.e.\ starting at $x=0$ (see Figure \ref{3regions}). In \eqref{SchBVP}, the Robin boundary condition (BC) at $x=0$ only fixes the ratio $\frac{\psi'(0)}{\psi(0)}$, hence an auxiliary Dirichlet (or Neumann) boundary value has to be invoked here. Its correct value will then be obtained by scaling the final equation using the inhomogeneous Robin BC at $x=1$. In contrast to \cite{ABAN11}, \eqref{SchBVP} includes the evanescent region (2), cf.\ Fig.\ \ref{3regions}, which still has to be formulated as a BVP (for stability reasons). This leads to the following domain decomposition and problem coupling for the auxiliary wave functions $\zeta,\,\chi,\,\varphi$:\\

\noindent
\underline{Step 1 -- IVP for $\zeta$ in region (1):}
\be \label{region1}
\left\{
\begin{array}{l}
\ds \eps^2 \zeta ''(x) + a(x) \zeta(x)=0\,, \quad x \in (0,x_c)\,, \\[3mm]
\ds \zeta(0) =1\,, \qquad\qquad\qquad \mbox{(auxiliary Dirichlet BC)}\\[3mm]
\ds \eps\zeta'(0) =-{i} \sqrt{a(0)} \,.\qquad \mbox{(due to the left BC in \eqref{SchBVP})}
\end{array}
\right.
\ee
\medskip

\noindent
\underline{Step 2 -- BVP for $\chi$ in region (2):}
\be \label{region2}
\left\{
\begin{array}{l}
\ds \eps^2 \chi ''(x) + a(x) \chi(x)=0\,, \quad x \in (x_c,x_d)\,, \\[3mm]
\ds \zeta'(x_c)\chi(x_c) - \zeta(x_c)\chi'(x_c)=0\,, \quad \mbox{(Robin BC for $\chi$ implies continuity of $\frac{\psi'}{\psi}$ at $x_c$)}\\[3mm]
\ds \eps\chi'(x_d) =1\,. \qquad\qquad\qquad\qquad\quad \mbox{(auxiliary Neumann BC)}
\end{array}
\right.
\ee
\medskip

\noindent
\underline{Step 3 -- IVP for $\varphi$ in region (3):}
\be \label{region3}
\left\{
\begin{array}{l}
\ds \eps^2 \varphi ''(x) + a(x) \varphi(x)=0\,, \quad x \in (x_d,1)\,, \\[3mm]
\ds \varphi(x_d) =\chi(x_d)\,, \qquad\:\: \mbox{(implies continuity of $\psi$ at $x_d$)}\\[3mm]
\ds \varphi'(x_d) =\chi'(x_d) \,.\qquad \mbox{(implies continuity of $\psi'$ at $x_d$)}
\end{array}
\right.
\ee
\medskip

\noindent
\underline{Step 4 -- scaling of the auxiliary wave functions:}
\be \label{psi-scaling}
\psi(x):=\left\{
\begin{array}{l}
\ds \beta\,\zeta(x)\,, \quad x \in (0,x_c)\,, \\[3mm]
\ds \alpha\,\chi(x)\,, \quad x \in (x_c,x_d)\,, \\[3mm]
\ds \alpha\,\varphi(x)\,, \quad x \in (x_d,1)\,,
\end{array}
\right.
\ee
with the scaling parameters $\alpha,\,\beta\in\CC$ still to be defined.

This procedure can be explained as follows: First we note that the BCs of \eqref{region1} imply $\eps\zeta'(0)+i\sqrt{a(0)}\zeta(0)=0$, just as in the left BC of the BVP \eqref{SchBVP}. Hence,
the IVP \eqref{region1} coincides on region (1) with the BVP \eqref{SchBVP}, except for the auxiliary Dirichlet BC $\zeta(0)=1$. The true solution of \eqref{SchBVP} satisfies instead $\psi(0)=\beta$ with some a-priori unknown $\beta\in\CC$. Hence, the auxiliary wave function $\zeta$ is related to $\psi$ by the scaling $\psi\big|_{[0,x_c]}=\beta\zeta$, as postulated in the first line of \eqref{psi-scaling}. Clearly, this implies $\frac{\psi'}\psi = \frac{\zeta'}\zeta$ on $[0,x_c]$. In the above Step 2, the Robin BC allows to carry over this relation to region (2): $\frac{\psi'}\psi = \frac{\chi'}\chi$, and the auxiliary wave function $\chi$ is related to $\psi$ by the scaling $\psi\big|_{[x_c,x_d]}=\alpha\chi$, with some $\alpha\in\CC$ to be determined. 
The initial conditions for the auxiliary wave function $\varphi$ in Step 3 imply $C^1$--continuity of $\psi$ when using again the scaling $\psi\big|_{[x_d,1]}=\alpha\varphi$.

So far, the wave function $\psi$ defined in  \eqref{psi-scaling} neither satisfies continuity at $x_c$ nor the right BC from \eqref{SchBVP}. Therefore we define the scaling parameters $\alpha,\,\beta\in\CC$ via
\begin{eqnarray}
  &&\qquad \alpha\,[\eps\varphi'(1)-i \sqrt{a(1)}\,\varphi(1)]=-2i \sqrt{a(1)} \,,\qquad \mbox{(due to the right BC in \eqref{SchBVP})}
  \label{alpha-eq}\\
  &&\qquad \beta\,\zeta(x_c)=\alpha\,\chi(x_c) \,.\qquad\qquad\qquad\qquad\qquad\:\:\: \mbox{(implies continuity of $\psi$ at $x_c$)}
  \label{beta-eq}
\end{eqnarray}

\begin{rem}\label{rem2.1}
The key aspect of the above algorithm is to prescribe in the BVP \eqref{region2} the continuity of $\frac{\zeta'}{\zeta}$ to $\frac{\chi'}{\chi}$ at $x_c$. Note that this continuity is invariant under the scaling \eqref{psi-scaling}. Hence it is inherited by $\frac{\psi'}{\psi}$, implying (with the continuity of $\psi$) the required $C^1$--continuity of $\psi$. The simpler alternative to prescribe in \eqref{region2} continuity of $\zeta$ to $\chi$ would typically be paired with a discontinuity of $\zeta'$ to $\chi'$ at $x=x_c$ (as a result of solving the BVP). Then, the scaling of \eqref{psi-scaling} would lead to an unwanted discontinuity of $\psi'$ at $x=x_c$.
\end{rem}

\medskip

\begin{lem}\label{BVP-3zones}
Let Hypothesis a3 be satisfied.
Then:
\begin{enumerate}
\item[(i)] The IVPs \eqref{region1}, resp. \eqref{region3} admit $\forall \eps >0$ unique solutions $\zeta\in W^{2,\infty}(0,x_c)\subset C^1[0,x_c]$, resp. $\varphi \in W^{2,\infty}(x_d,1)\subset C^1[x_d,1]$.
\item[(ii)] The BVP \eqref{region2} has $\forall \eps >0$ a unique solution $\chi\in W^{2,\infty}(x_c,x_d)\subset C^1[x_c,x_d]$.
\item[(iii)] The scaling parameters $\alpha,\,\beta\in\CC\setminus\{0\}$ are well-defined by \eqref{alpha-eq}, \eqref{beta-eq}.
\end{enumerate}
\end{lem}
\medskip

\debproof Part (i) is straightforward. For (ii),  let us first consider the IVP \eqref{region1}. Its unique solution $\zeta$ has the property: The values $\zeta(x_c)$ and $\zeta'(x_c)$ are linearly independent over the field $\RR$. Otherwise, the backward IVP (starting at $x_c$) would yield ``final values'' $\zeta(0)$ and $\zeta'(0)$ that are linearly dependent over $\RR$, which is in contradiction with the initial condition in \eqref{region1}.

To solve the BVP \eqref{region2}, let $\chi_1$, $\chi_2$ be a (real valued) basis of solutions for that Schr\"odinger equation on $(x_c,x_d)$, with
\begin{eqnarray*}
  &&\chi_1(x_c)=1\,,\qquad\chi_1'(x_c)=0\,,\\
  &&\chi_2(x_c)=0\,,\qquad\chi_2'(x_c)=1\,.
\end{eqnarray*}
Setting $\chi=c_1\chi_1+c_2\chi_2$ with some $c_1,\,c_2\in\CC$, the BCs of \eqref{region2} give rise to the following linear equation:
\begin{equation}\label{chi-coefficients}
  \left(
\begin{array}{cc}
\zeta'(x_c)&-\zeta(x_c)\\
\eps\chi_1'(x_d)&\eps\chi_2'(x_d)
\end{array}
\right)\:
  \left(
\begin{array}{c}
c_1\\
c_2
\end{array}
\right)\:=\:
  \left(
\begin{array}{c}
0\\
1
\end{array}
\right)\,.
\end{equation}
The determinant of this system satisfies 
$\eps[\zeta'(x_c)\chi_2'(x_d)+\zeta(x_c)\chi_1'(x_d)]\ne0$, since $\chi_{1,2}'(x_d)$ $\in\RR$, 
but $\zeta(x_c)$ and $\zeta'(x_c)$ are linearly independent over $\RR$. Hence, \eqref{chi-coefficients} is uniquely solvable for $c_1,\,c_2$, thus providing the unique solution to \eqref{region2}.\\

\noindent
For part (iii) we shall first argue that \eqref{alpha-eq} yields a well-defined $\alpha\in \CC\setminus\{0\}$. To this end we shall prove that $\eps\varphi'(1)\ne i \sqrt{a(1)}\,\varphi(1)$, using the quantum mechanical current of the model \eqref{SchBVP}: 
\begin{equation}\label{current}
  j(x):=\eps \Im [\bar\psi(x)\,\psi'(x)]\,.
\end{equation}
Assume now that $\eps\varphi'(1)= i \sqrt{a(1)}\,\varphi(1)$. Then, \eqref{psi-scaling} implies on the one hand
\begin{equation}\label{current1}
  j(1)=\eps |\alpha|^2 \,\Im [\bar\varphi(1)\,\varphi'(1)]
  =|\alpha|^2 \,\sqrt{a(1)}\,|\varphi(1)|^2\ge0\,.
\end{equation}
But, on the other hand, \eqref{region1} yields
$$
  j(0)=\eps |\beta|^2 \Im [\bar\zeta(0)\,\zeta'(0)]
  =-|\beta|^2 \sqrt{a(0)}\,|\zeta(0)|^2\le0\,.
$$
Since the current in a stationary quantum model is constant in $x$, this implies $j\equiv 0$.
Since $a(1)>0$ and $\alpha\ne0$ (otherwise $\psi(1)=\psi'(1)=0$ would contradict the BC at $x=1$ in \eqref{SchBVP}), \eqref{current1} shows that $\varphi(1)=0$, and hence --by our assumption-- $\varphi'(1)=0$. But this leads to a contradiction in the BC at $x=1$ in \eqref{SchBVP}. Hence, \eqref{alpha-eq} yields indeed a unique $\alpha\in \CC\setminus\{0\}$.

Finally, \eqref{beta-eq} yields a well-defined $\beta\in\CC$ since $\zeta(x_c)$ and $\zeta'(x_c)$ are linearly independent over $\RR$ (as shown in part (i) above). Moreover, $\beta\ne0$ since $\chi(x_c)=c_1\ne0$ (otherwise the first line of \eqref{chi-coefficients} would also yield $c_2=0$).
\finproof

\medskip
The above lemma shows that the \emph{domain decomposition algorithm} \eqref{region1}-\eqref{beta-eq} yields a unique function $\psi$ that is piecewise in $W^{2,\infty}$ and piecewise (on the three regions) a solution to the Schr\"odinger equation \eqref{SchBVP}. Moreover, one has the proposition:
\vspace{0.2cm}
\debprop\label{prop2.5}
Let Hypothesis a3 be satisfied.
Then the function $\psi$ obtained from \eqref{region1}-\eqref{beta-eq} belongs to $W^{2,\infty}(0,1)$ and is the unique solution of \eqref{SchBVP} (as guaranteed by Proposition \ref{Prop1.1}).
\finprop\\
\vspace{0.2cm}

\debproof
  The matching conditions in \eqref{region2} and \eqref{beta-eq} imply $C^1$--continuity of $\psi$ at $x_c$, and the initial conditions in \eqref{region3} imply $C^1$--continuity of $\psi$ at $x_d$. Hence, $\psi \in C^1[0,1]$, and this proves the claim.
\finproof

\Section{Numerical analysis of the hybrid WKB-method} \label{SEC3}
To keep the presentation simple we shall consider here only the two-zone model of \S\ref{SEC2.1}. It has a coefficient function $a(x)$ that corresponds to Figure \ref{2regions}. In the following subsections we shall thus study step by step the numerical errors obtained when solving the BVP \eqref{region1_ev} in the evanescent region with a multiscale WKB-FEM and the IVP \eqref{region2_ev} in the oscillatory region with the marching method introduced in \cite{ABAN11}. We shall always assume that the phase function $\int^x\sqrt{|a(y)|}\,dy$ in the WKB-functions \eqref{osWKB-fct}, \eqref{evWKB-fct} can be computed exactly, e.g., this holds for piecewise linear $a(x)$. Otherwise, an additional quadrature error of the phase would need to be included in our subsequent analysis.

\subsection{Variational formulation  for the evanescent region BVP \eqref{region1_ev}} \label{SEC311}
Let us introduce in this section the variational formulation of the evanescent region problem (\ref{region1_ev}) and study the well-posedness of the problem. As pointed out previously, we consider in the current paper situations with an abrupt potential jump, avoiding turning points, such that we shall suppose:\\

{\bf Hypothesis A} {\it Let $V \in W^{2,\infty}(0,x_d)$ and $E>0$ satisfy the bounds
$$
0<\tau_{ev}  \le -a(x):=V(x)-E\le M_{ev}\,, \quad \forall x \in (0,x_d)\,.
$$
Furthermore let in the following $0 < \eps <1$ be arbitrary.}\\

We are now searching for a weak solution of (\ref{region1_ev}) in the Hilbert space 

$$
{\mathcal{V}}:=H^1(0,x_d)\,, \quad (\chi,\theta)_{\mathcal V}:=(\chi,\theta)_{L^2(0,x_d)}+\eps^2(\chi',\theta')_{L^2(0,x_d)}\,.
$$
This $\eps$--dependent norm gives rise to the following weighted Sobolev embedding, where the Gagliardo-Nirenberg inequality for bounded domains is used in the first estimate:
\be\label{Sobolev-emb}
\begin{array}{ccc}
 \ds \eps \|\chi\|_{C[0,x_d]}^2 &\le&\ds C \|\chi\|_{L^2(0,x_d)}\,\big(\eps\|\chi'\|_{L^2(0,x_d)}\big)+C \eps \|\chi\|^2_{L^2(0,x_d)}\\[3mm]
&  \le&\ds C \big(\|\chi\|_{L^2(0,x_d)}^2 + \eps^2\|\chi'\|_{L^2(0,x_d)}^2\big) =C \|\chi\|^2_{\mathcal V}\,.
\end{array}
\ee

The variational formulation reads: Find $\chi \in {\mathcal{V}}$, solution of
\be \label{VF_ev}
b(\chi,\theta)=L(\theta)\,, \quad \forall \theta \in \mathcal{V}\,,
\ee
with the sesquilinear form $b:\mathcal{V} \times \mathcal{V} \rightarrow \CC$ and the antilinear form $L:\mathcal{V} \rightarrow \CC$ defined as
\be \label{bL_ev}
\begin{array}{lll}
\ds b(\chi,\theta)\!\!\!\!&:=&\!\!\!\! \ds \eps^2 \int_0^{x_d}\!\! \chi'(x)\, \overline{\theta}'(x)\, dx -\int_0^{x_d}\!\! a(x)\, \chi(x)\, \overline{\theta}(x)\, dx + \eps \sqrt{|a(0)|}\, \chi(0)\,\overline{\theta}(0)\,, \quad \forall \chi,\theta \in \mathcal{V}\\[3mm]
\ds L(\theta)\!\!\!\!&:=&\!\!\!\! \ds \eps \,  \overline{\theta}(x_d)\,, \quad \forall \theta \in \mathcal{V}\,.
\end{array}
\ee
The BVP \eqref{region1_ev} is a standard elliptic problem, meaning that the forms $b(\cdot,\cdot)$ and $L(\cdot)$ are continuous and $b(\cdot,\cdot)$ is coercive, {\it i.e.} there exists a constant $C>0$ independent of $\eps$, such that for all $\chi,\theta \in \mathcal{V}$ one has
$$
|b(\chi,\theta)|\le C\, \|\chi\|_\mathcal{V}\, \|\theta\|_\mathcal{V}\,, \quad |L(\theta)| \le C\sqrt\eps
\, \|\theta\|_\mathcal{V}\,, \quad b(\theta,\theta) \ge \min \{1, \tau_{ev} \}\, \|\theta\|_\mathcal{V}^2\,.
$$
The Lax-Milgram theorem implies then for each $\eps >0$ the existence and uniqueness of a weak solution $\chi \in \mathcal{V}$ of \eqref{VF_ev}. We have moreover the following lemma:\\

\begin{lemma}\label{lem3.1}
Let Hypothesis A be satisfied. Then the weak solution $\chi \in {\mathcal V}$ of \eqref{VF_ev} or \eqref{region1_ev} belongs even to $\in H^2(0,x_d) \hookrightarrow C^1[0,x_d]$ and satisfies the following estimates, with a constant $C>0$ independent on $\eps$
\be \label{CHI1}
\|\chi\|_{L^2(0,x_d)}^2 \le C \eps\,, \quad \eps^2\, \|\chi'\|_{L^2(0,x_d)}^2 \le C \eps\,,
\ee
as well as
\be \label{CHI2}
\|\chi\|_{C[0,x_d]} \le C\,, \quad \eps\, \|\chi'\|_{C[0,x_d]} \le C\,.
\ee
\end{lemma}
\debproof
%
The Lax-Milgram theorem yields immediately
$$
  \|\chi\|_{\mathcal V} \le \frac{C\sqrt\eps}{\min\{1,\tau_{ev}\}}\,,
$$
which implies \eqref{CHI1} and, with \eqref{Sobolev-emb}, the first inequality of \eqref{CHI2}.

To show the second estimate of \eqref{CHI2}, we observe that
$$
\begin{array}{lll}
\eps^2 (\chi'(x))^2&=&\ds \eps^2 (\chi'(0))^2 + 2 \int_0^x \chi'(y)\, [\eps^2\, \chi^{''}(y)]\, dy\\[3mm]
&=& \ds |a(0)|\, |\chi(0)|^2 - 2 \int_0^x a(y) \chi'(y)\, \chi(y)\, dy\\[3mm]
&\le& \ds C+C\|\chi'\|_{L^2}\, \|\chi\|_{L^2} \le C+C\eps \|\chi'\|_{L^2}^2+{C \over \eps}\, \|\chi\|^2_{L^2}\le C\,,
\end{array}
$$
where we used the other, just proved estimates.
\finproof
\subsection{Convergence analysis for the WKB-FEM method in the evanescent region} \label{SEC32}

The multi-scale WKB-FEM method we shall use for an efficient resolution of the evanescent region problem (\ref{region1_ev}) is based on a specific choice of WKB-basis functions from \eqref{evWKB-fct}. In more detail, the Hilbert space $\mathcal{V}$ will be approximated by an appropriate finite-dimensional Hilbert space $\mathcal{V}_h \subset {\mathcal V}$, spanned by well chosen basis functions, and the continuous problem \eqref{VF_ev} will be approximated by the following discrete problem: Find $\chi_h \in {\mathcal{V}_h}$, solution of
\be \label{VF_ev_h}
b(\chi_h,\theta_h)=L(\theta_h)\,, \quad \forall \theta_h \in \mathcal{V}_h\,.
\ee
To introduce the finite-dimensional space $\mathcal{V}_h$, let us partition the interval $[0,x_d]$ into $0= x_1 < x_2 <  \cdots < x_N = x_d$ and denote the mesh size by $h_n:=x_{n+1}-x_n$ as well as $h:= \max_{n=1, \cdots,N-1} \{h_n\}$.
The appropriate Hilbert space $\mathcal{V}_h$ is then defined as
$$
\mathcal{V}_h:= \left\{ \theta_h \in \mathcal{V} \,\, \Big| \,\, \theta_h(x)= \sum_{n=1}^N z_n\, \zeta_n(x)\,, \quad z_n \in \CC \right\}\,,
$$
with the WKB-hat functions defined as
\be \label{hat}
\zeta_n(x):=
\left\{
\begin{array}{ll}
v_{n-1}(x)\,,& x \in [x_{n-1},x_n]\,,\\[3mm]
w_{n}(x)\,,& x \in [x_{n},x_{n+1}]\,,\\[3mm]
0\,,&  \mbox{otherwise}\,.
\end{array}
\right.
\ee
Here we used the notation
\be \label{NOTZ}
\begin{array}{lll}
\ds w_n(x):= \alpha_n(x)\, q_{n}(x)&;&\ds v_n(x):= \beta_n(x)\, q_{n+1}(x)\,,\\[3mm]
\ds \alpha_n(x):=- {\sinh \sigma_{n+1}(x) \over \sinh \gamma_n}&;& \ds \beta_n(x):= {\sinh \sigma_{n}(x) \over \sinh \gamma_n}\,, \\[3mm]
\ds \sigma_n(x):= { 1 \over \eps} \int_{x_n}^x \sqrt{|a(y)|}\, dy&;&\ds \gamma_n:= { 1 \over \eps} \int_{x_n}^{x_{n+1}} \sqrt{|a(y)|}\,dy\,,\\[3mm]
\ds q_n(x):= \frac{(V(x_n)-E)^{1/4}}{(V(x)-E)^{1/4}}\,.
\end{array}
\ee
Assuming Hypothesis A, $\zeta_n$ is piecewise in $W^{2,\infty}(x_j,x_{j+1})$ $\forall j\in\{1,...,N-1\}$ and globally in $W^{1,\infty}(0,x_d)\hookrightarrow C[0,x_d]$. Note that both components $v_n$ and $w_n$ of these (non-standard) hat functions are linear combinations of the evanescent WKB-functions of first order, i.e. $\varphi_1^{ev}$ given in \eqref{evWKB-fct}. Hence these hat functions are solutions of our Schr\"odinger equation up to an error of order $\mathcal{O}(\eps^2)$, {\it i.e.}
\begin{eqnarray*}
\eps^2 \zeta_n^{''}(x) + a(x)\, \zeta_n(x) = \eps^2 \left[ {5 \over 16}\, {(V'(x))^2 \over (V(x)-E)^2 } + {1 \over 4} {V^{''}(x) \over E-V(x)}\right]\, \zeta_n(x)\,, &&\\
\forall x \in (x_{n-1},x_n)\cup (x_n,x_{n+1})\,.&&
\end{eqnarray*}
This peculiarity signifies that the hat functions incorporate already some essential information about the  solutions we are searching for, leading to a scheme which will be asymptotically correct in the limit $\eps \rightarrow 0$, as will be seen later on.\\
For later purposes let us introduce here the differential operator
$$
{\mathcal A}_{\eps}(\xi) :=  -\eps^2 \xi^{''}(x) - a(x)\, \xi(x) +\eps^2 r(x)\, \xi(x)\,,
$$
with the function
$$
r(x) :=  \ds  {5 \over 16}\, {(V'(x))^2 \over (V(x)-E)^2 } + {1 \over 4} {V^{''}(x) \over E-V(x)}\,.
$$
For each $\xi \in {\mathcal V}_h$ one has ${\mathcal A}_{\eps}(\xi) =0$ in every interval $I_n:=(x_n,x_{n+1})$. But we note that ${\mathcal A}_{\eps}$ cannot be applied globally on $(0,x_d)$, as functions in ${\mathcal V}_h$ typically have discontinuous derivatives at the grid points $x_n$.\\
We also remark that, due to Lax-Milgram's theorem, the discrete problem \eqref{VF_ev_h} admits $\forall \eps >0$ also a unique solution $\chi_h \in {\mathcal V}_h$.\\

The aim is now to investigate the error between the exact solution of \eqref{VF_ev}, denoted by $\chi_{ex}$, and the solution of the discrete problem \eqref{VF_ev_h}, denoted by $\chi_h$. We denote by $\Pi_h^{\eps} \chi_{ex} \in \mathcal{V}_h$ the interpolant of the exact solution in the finite dimensional Hilbert space $\mathcal{V}_h$, {\it i.e.}
\be\label{projection}
\Pi_h^{\eps} \chi_{ex} (x) := \sum_{n=1}^N \chi_{ex}(x_n)\, \zeta_n(x)\,, \quad \forall x \in [0,x_d]\,.
\ee
Then, the numerical error can be split as follows
$$
e_h(x):= \chi_{ex}(x)-\chi_h(x) = \left(\chi_{ex}(x)- \Pi_h^{\eps} \chi_{ex}\right) + \left(\Pi_h^{\eps} \chi_{ex}-\chi_h \right)=:e_h^1(x)+e_h^2(x)\,,
$$
where $e_h^1$ corresponds to the interpolation error (consistency) and $e_h^2$ is the stability error. These two error parts shall be now estimated separately.
\subsubsection{Consistency error estimate} \label{SEC321}
The goal of this section is to estimate the interpolation error $e_h^1(x):=\chi_{ex}(x)-\Pi_h^{\varepsilon} \chi_{ex}(x)$ in the $\mathcal{V}$-norm.\\

To this end, note that the equation satisfied by $e_h^1$ in $I_n$ is 
$$
\left\{
\begin{array}{l}
{\mathcal A}_{\eps}(e_h^1) = \eps^2\, r(x)\, \chi_{ex}(x)\,, \qquad \forall x \in I_n\,,\\[3mm]
\ds e_h^1(x_n)=e_h^1(x_{n+1})=0\,.
\end{array}
\right.
$$
The variation of constants method, i.e.\ writing $e_h^1(x)=c_1(x)\, w_n(x)+c_2(x)\, v_n (x)$ in $I_n$ leads after some lengthy but straightforward computations (see \cite{Cla} for the oscillatory case) to the following explicit expressions for the error function 
\be \label{CAP3err}
e_h^1(x)=\mathcal{E}_1(x)+\mathcal{E}_2(x)\, , \qquad x \in I_n\,,
\ee
with
$$
\begin{array}{lll}
\ds\mathcal{E}_1(x)&\!\!\! =&\!\!\! \ds -\frac{\varepsilon}{(V(x)-E)^{1/4}}\, { \sinh \sigma_{n+1} (x)\over \sinh \gamma_n  }\, \int_{x_n}^x \frac{r(y)\, \chi_{ex}(y)}{(V(y)-E)^{1/4}} \sinh \sigma_n(y)\, dy\, ,\\[5mm]
\ds\mathcal{E}_2(x)&\!\!\! =&\!\!\! \ds  \frac{\varepsilon}{(V(x)-E)^{1/4}} \, \frac{\sinh \sigma_n(x)}{\sinh \gamma_n} \, \int_{x_{n+1}}^{x} \!\! \frac{r(y)\, \chi_{ex}(y)}{(V(y)-E)^{1/4}} \sinh \sigma_{n+1}(y)\, dy\, .
\end{array}
$$
Differentiating the interpolation error function yields
\be \label{CAP3derr}
(e_{h}^1)'(x)=\mathcal{D}_1(x)+\mathcal{D}_2(x)+\mathcal{D}_3(x)+\mathcal{D}_4(x)
= \mathcal{D}_2(x)+\mathcal{D}_4(x)\, ,
\ee
with
$$
\begin{array}{lll}
\ds \mathcal{D}_1(x)&\!\!\! =&\!\!\! -\mathcal{D}_3(x) =\ds -\eps\, { \sinh \sigma_{n+1} (x)\over \sinh \gamma_n  }\, \frac{\sinh \sigma_n(x)}{\sqrt{V(x)-E}}\, r(x)\, \chi_{ex}(x)\,,\\[5mm]
\ds \mathcal{D}_2(x)&\!\!\!=&\!\!\! \ds - {\eps \over \sinh \gamma_n} \, \left[ {(V(x)-E)^{1/4} \over \eps } \cosh \sigma_{n+1}(x) - {V'(x)\, \sinh \sigma_{n+1}(x) \over 4\, (V(y)-E)^{5/4} } \right]  \\[5mm]
&& \times \,\ds\int_{x_n}^x \frac{r(y) \, \chi_{ex}(y)\, \sinh \sigma_n(y)}{(V(y)-E)^{1/4}} \, dy\, , \\[5mm]
\ds \mathcal{D}_4(x)&\!\!\! =&\!\!\! \ds  {\eps \over \sinh \gamma_n} \, \left[ {(V(x)-E)^{1/4} \over \eps } \cosh \sigma_{n}(x) - {V'(x)\, \sinh \sigma_{n}(x) \over 4\, (V(y)-E)^{5/4} } \right]  \\[5mm]
&& \times\,\ds \int_{x_{n+1}}^x \frac{r(y) \, \chi_{ex}(y)\, \sinh \sigma_{n+1}(y)}{(V(y)-E)^{1/4}} \, dy\, .
\end{array}
$$
In order to estimate the interpolation error in the $\mathcal{V}$-norm, we shall investigate each of these terms separately. In this study, the behaviour of the following functions is very important:
$$
\Theta_{ss}(x):=\frac{\sinh (\sigma_n(x))\,\sinh (-\sigma_{n+1}(x)) }{\sinh \gamma_n}\,, \quad \Theta_{sc}(x):=\frac{\sinh (\sigma_n(x))\,\cosh (-\sigma_{n+1}(x)) }{\sinh \gamma_n}\,,
$$
$$
\Theta_{cs}(x):=\frac{\cosh (\sigma_n(x))\,\sinh (-\sigma_{n+1}(x)) }{\sinh \gamma_n}\,, \qquad \forall x \in I_n\,. 
$$
Next we shall use $\left|\frac{r(y)}{(V(y)-E)^{1/4}}\right|\le C$ and the fact that $\sinh \sigma_n(x)$ and $\sinh(- \sigma_{n+1}(x))$ are both non-negative on $I_n$ and, respectively, increasing and decreasing. Then one can show for $x \in I_n$:
$$
|\mathcal{E}_1(x)| \le C \eps \, \Theta_{ss}(x)\,  \int_{x_n}^x |\chi_{ex}(y)|\, dy\,, \qquad |\mathcal{E}_2(x)| \le C \eps \, \Theta_{ss}(x)\,  \int_{x}^{x_{n+1}} |\chi_{ex}(y)|\, dy \,,
$$
and
$$
|\mathcal{D}_2(x)|\le C\, \Theta_{sc}(x) \int_{x_n}^x |\chi_{ex}(y)|\, dy + C \eps\, \Theta_{ss}(x) \int_{x_n}^x |\chi_{ex}(y)|\, dy \,,
$$
$$
|\mathcal{D}_4(x)|\le C\, \Theta_{cs}(x) \int_{x}^{x_{n+1}} |\chi_{ex}(y)|\, dy + C \eps\, \Theta_{ss}(x) \int_{x}^{x_{n+1}} |\chi_{ex}(y)|\, dy \,.
$$
In the above estimates, the constant $C$ depends only on our data $a(x)$ and $E$, but not on $\eps$ and $\chi_{ex}$. Using $\sigma_n(x)-\sigma_{n+1}(x)=\gamma_n$ we easily find
$$
0 \le \Theta_{ss}(x) \le \frac{\cosh \gamma_n -1 }{2\, \sinh \gamma_n} \le {1 \over 2}\,, \quad \forall x \in I_n\,, 
$$
$$
0 \le \Theta_{sc}(x) = {1 \over 2} +\frac{\sinh (\sigma_n(x) + \sigma_{n+1}(x)) }{2\,  \sinh \gamma_n} \le 1\,,\quad \forall x \in I_n\,, 
$$
$$
0 \le \Theta_{cs}(x) = {1 \over 2}-\frac{\sinh (\sigma_n(x) +\sigma_{n+1}(x)) }{2\,  \sinh \gamma_n}\le 1\,, \quad \forall x \in I_n\,.
$$
With the asymptotic behaviour
$$
\frac{\cosh \xi -1}{2\, \sinh \xi} \;\stackrel{\xi \sim 0}{\sim} \;
{\xi \over 4}\,, \quad \frac{\cosh \xi -1}{2\, \sinh \xi}  \stackrel{\xi \rightarrow \pm \infty}{\longrightarrow} \pm {1 \over 2}
$$
we obtain (using $|\gamma_n|\le C\frac{h}{\eps}$)
$$
  \eps | \Theta_{ss}(x)| \le C \min \{\eps,h\}\,.
$$
With Lemma \ref{lem3.1} this permits to prove the following lemma:\\

\begin{lemma}\label{th:e1-est}
Let Hypothesis A be satisfied. Then the following estimates hold for the interpolation error $e_h^1\in {\mathcal V}_h \subset C[0,x_d]$ of the exact solution $\chi_{ex} \in \mathcal{V}$ of \eqref{VF_ev}:
$$
\|e_h^1\|_{L^2(0,x_d)} \le C \sqrt{\eps}\, h\, \min \{ \eps,h\}\,, \quad \eps\,\|(e_h^1)'\|_{L^2(0,x_d)} \le C \eps^{3/2}\, h\,,
$$
$$
\|e_h^1\|_{C[0,x_d]} \le C \,\sqrt{h}\, \min \{ \eps^{3/2},h^{3/2}\}\,, \quad \eps  \|(e_h^1)'\|_{L^{\infty}(0,x_d)}\le C \eps\, \sqrt{h}\,\min \{ \sqrt{\eps},\sqrt{h}\} \,.
$$
\end{lemma}
\subsubsection{Stability error estimate} \label{SEC322}
Let us now come to the ${\mathcal V}$-estimates of the stability error
$$
e_h^2(x):=\Pi_h^{\eps} \chi_{ex}-\chi_h \,.
$$
For this study, we remark that
$$
b(\chi_{ex}-\chi_h,\theta_h)=0 \qquad \forall \theta_h \in {\mathcal V}_h\,,
$$
implying with the choice $\theta_h= e_h^2 \in {\mathcal V}_h$ that
\be \label{ERR1}
b(e_h^2,e_h^2)=b(\Pi_h^{\eps} \chi_{ex}-\chi_h,e_h^2)=b(\Pi_h^{\eps} \chi_{ex}-\chi_{ex},e_h^2)=-b(e_h^1,e_h^2)\,.
\ee
Using now the coercivity and the boundedness of the sesquilinear form $b$, yields with a constant $C>0$ independent of $\eps$, that
%
$$
\min \{1, \tau_{ev} \}\, \|e_h^2\|_{\mathcal V} \le
b(e_h^2,e_h^2) \big/ \|e_h^2\|_{\mathcal V} 
\le C \,  \|e_h^1\|_{\mathcal V}\,,
$$
but this estimate can be further improved for $0<\eps \ll 1$. For this, let us study the right hand side of \eqref{ERR1} in more detail. As $e_h^2 \in {\mathcal V}_h$ we have that ${\mathcal A}_{\eps}(e_h^2)=0$ on every interval $I_n$, implying thus
$$
-b(e_h^1,e_h^2)=\eps^2 \sum_{n=1}^{N-1} \int_{x_n}^{x_{n+1}} e_h^1\, (\overline{e_h^2})^{''}\, dx + \int_0^{x_d} a(x)\, e_h^1\, \overline{e_h^2}\, dx= - \eps^2 \int_0^{x_d} r(x) e_h^1\, \overline{e_h^2}\, dx\,,
$$
were we used $e_h^1(x_n)=0$ for the integration by parts.
Thus one obtains from \eqref{bL_ev}, \eqref{ERR1}:
$$
\eps^2 \|(e_h^2)'\|^2_{L^2} + \|e_h^2\|^2_{L^2} + \eps |e_h^2(0)|^2 \le C\, b(e_h^2,e_h^2)=-C\, b(e_h^1,e_h^2) \le C \eps^2\,\|e_h^1\|_{L^2}\, \|e_h^2\|_{L^2}\,, 
$$
and in particular $\|e_h^2\|_{L^2}\le C \eps^2\,\|e_h^1\|_{L^2}$.
Using also the Sobolev embedding \eqref{Sobolev-emb},
this implies
\begin{eqnarray}\label{e2-est}
\|e_h^2\|^2_{L^2(0,x_d)} + \eps^2 \|(e_h^2)'\|^2_{L^2(0,x_d)}  + \eps \|e_h^2\|^2_{C[0,x_d]} &\le& C \eps^4 \|e_h^1\|_{L^2(0,x_d)}^2\, \nonumber\\
&\le& C\eps^5 h^2 \min\{\eps^2,h^2\}\,.
\end{eqnarray}

Finally we shall now derive an $L^\infty$--bound on $(e_h^2)'$, using the bound on $\|e_h^2\|_{C[0,x_d]}$: 
Since $e_h^2 \in\mathcal V_h$, we have 
$$
  e_h^2(x) = \sum_{n=1}^N e_h^2(x_n) \,\zeta_n(x)\,,\qquad
  (e_h^2)' (x) = \sum_{n=1}^N e_h^2(x_n) \,\zeta_n'(x)\, \quad \mbox{a.e. on }[0,x_d]\,.
$$
To bound $\|\zeta_n'\|_ {L^{\infty}}$, one refers to \eqref{NOTZ} and observes that  $0\le\alpha_n(x)\le1$ as well as
$$
  0\le-\alpha_n'(x)\le \frac{\sqrt{M_{ev}}}{\eps} \coth \gamma_n \le C\big(\frac1\eps+\frac1h\big)   \quad\mbox{on }I_n\,.
$$
And analogous estimates hold for $\beta_n(x)$. From \eqref{e2-est} we then obtain
\be\label{e2-est2}
  \eps \|(e_h^2)'\|_{L^{\infty}(0,x_d)}\le C \eps \|e_h^2\|_{C[0,x_d]}\big(\frac1\eps+\frac1h\big) \le C \eps^3\, h \,.
\ee

\subsubsection{Convergence results for the WKB-FEM method} \label{SEC323}
To summarize, we shall put together both error contributions (from Lemma \ref{th:e1-est} and from the stability estimates \eqref{e2-est}, \eqref{e2-est2}) in the following theorem. It turns out that the consistency error is dominant here:\\

\begin{theorem} \label{EV_conv} {\bf (Convergence WKB-FEM)}
Let Hypothesis A be satisfied. Then the following estimates hold for the numerical error between the exact solution $\chi_{ex} \in \mathcal{V}$ of \eqref{VF_ev} and the numerical solution $\chi_h \in {\mathcal V}_h$ of \eqref{VF_ev_h}:
\bean
  \|e_h\|_{L^2(0,x_d)} &\le& C \sqrt{\eps}\, h\, \min \{ \eps,h\}\,, \qquad\quad \eps\,\|e_h'\|_{L^2(0,x_d)} \le C \eps^{3/2}\, h\,,\\
  \|e_h\|_{C[0,x_d]} &\le& C \,\sqrt{h}\, \min \{ \eps^{3/2},h^{3/2}\}\,,\quad \eps  \|e_h'\|_{L^{\infty}(0,x_d)}\le C \eps\, \sqrt{h}\,\min \{ \sqrt{\eps},\sqrt{h}\} \,.
\eean
\end{theorem}

\subsection{Vectorial IVP for the oscillatory region} \label{SEC33}
In this subsection we shall first investigate the IVP \eqref{region2_ev} on the continuous level, on the interval $(x_d,1)$. Following \cite{ABAN11}, we shall rewrite \eqref{region2_ev} in vectorial form, done via a non-standard transformation that is appropriate for the numerical WKB-marching method. For the subsequent analysis let us make the following assumptions on the potential:
\vspace{0.2cm}

\noindent {\bf Hypothesis B} {\it Let $V \in C^5[x_d,1]$ and $E>0$ satisfy the bounds
$$
0<\tau_{os}  \le a(x):=E-V(x)\le M_{os}\,, \quad \forall x \in [x_d,1]\,.
$$
Moreover let $0 < \eps \le \eps_0$ be arbitrary, with some $\eps_0$ such that
$$
  0<\eps_0<\eps_1:=\min\left\{1,\,\min_{x_d\le x\le1} [a(x)^{1/4}\beta_+(x)^{-1/2}] \right\}.
$$
}\\
\vspace{0.2cm}

In this definition, $\beta_+$ denotes the non-negative part of $\beta$. Hence $\beta_+(x)^{-1/2}$ may take the value $\infty$. We note that the above restriction on $\eps$ guarantees that the phase function of $\varphi_2^{os}(x)$ (cf. \eqref{osWKB-fct}) is strictly increasing. Moreover, the resulting positivity of the function $\sqrt a-\eps^2\beta$ will be crucial for the WKB-marching method in \S\ref{SEC331}.
\\

Following \cite{ABAN11} it is convenient to pass from the second-order differential equation to a system of first-order, introducing the following vector notation for the wave function $\varphi(x)$ on $[x_d,1]$:
\be\label{def-U}
  \quad U(x)=\left( \begin{array}{c}
u_1\\
u_2
\end{array}
\right) :=  \left( \begin{array}{c}
a^{1/4}\varphi(x)\\
\frac{\eps (a^{1/4}\varphi)'(x)}{\sqrt{a(x)}}
\end{array}
\right) =  \left( \begin{array}{c}
a^{1/4}\varphi(x)\\
\eps \big(\frac14 a^{-5/4} a' \varphi + a^{-1/4}\varphi'\big)(x)
\end{array}
\right)\,.
\ee
The norm of $U$ is equivalent to the norm of the vector $(\varphi,\,\eps\varphi')^\top$. Indeed, the transformation matrix between these two vectors reads
\be \label{TRANSFO}
 \quad  A(x):=\left( \begin{array}{cc}
a^{1/4}(x)&0\\
\frac{\eps}{4} a^{-5/4}(x)a'(x) & a^{-1/4}(x)
\end{array}
\right) \quad \textrm{i.e.} \quad U(x)=A(x) \left(\begin{array}{c}\varphi\\\eps\varphi' \end{array}\right)\,,
\ee
where the matrix $A$ and its inverse are bounded, uniformly w.r.t.\ $x$ and $\eps$, due to Hypothesis B.\\

Let $\varphi_{ex}\in W^{2,\infty}(x_d,1)$ be the exact solution of \eqref{region2_ev} as guaranteed by Lemma \ref{PB-2zones}. In the above vector notation it will be denoted by $U_{ex}(x)$ or simply $U(x)$, and is solution to the system 
\begin{equation} \label{EQU}
\left\{
\begin{array}{l}
\ds U'(x)=\left[ \frac{1}{\eps} A_0(x)+\eps A_1(x)\right] U(x)\,,\quad x_d<x<1\,, \\[3mm]
\ds U(x_d)=A(x_d+)\, (\chi_{ex}(x_d), 1)^\top\,,
\end{array}
\right.
\end{equation}
with the two matrices
$$
A_0(x):=
\sqrt{a(x)}
\left(
\begin{array}{cc}
0&1\\
-1&0
\end{array}
\right)
\,; \quad A_1(x):=
\left(
\begin{array}{cc}
0&0\\
2 \beta(x)&0
\end{array}
\right)\, .
$$
Here, $\beta =-\frac{1}{2a^{1/4}}(a^{-1/4})''$ which was already defined in \eqref{WKB-fct}, and the matrix element $a(x_d+)$ of $A(x_d+)$ denotes the right-sided limit of $a$ at the jump discontinuity $x_d$. We also use the analogous notation for $a'(x_d+)$.

In the sequel we shall need an a-priori estimate on this solution. The upper bound was already given in \S2.1 of \cite{ABAN11}. But for the scaling Step 3 we shall also need an $\eps$-uniform lower bound on the solution:
\vspace{0.2cm}
\begin{lemma}\label{U-estimate}
Let Hypothesis B hold. Then, the ODE \eqref{EQU} admits a unique solution $U \in W^{1,\infty}(x_d,1)$, which satisfies
\be\label{U-est}
  \|U(x_d)\| \exp\Big[ -\eps \int_{x_d}^x  |\beta(y)| dy\Big] \le \|U(x)\| \le 
  \|U(x_d)\| \exp\Big[ \eps \int_{x_d}^x  |\beta(y)| dy\Big]\,,\quad x_d\le x\le1\,.
\ee
Thus, there exist constants $C_3,C_4>0$ independent on $\eps$ such that
\be\label{Uex-bound}
  C_3\le \|U\|_{C[x_d,1]} \le C_4\,,\qquad \forall 0<\eps\le\eps_0\,.
\ee
\end{lemma}\\
\vspace{0.2cm}

\debproof
For the norm 
$\|U\|^2 := |u_1|^2 + |u_2|^2$ we compute for \eqref{EQU}:
$$
  \left|\frac{d}{dx} \|U(x)\|^2\right| = 
  \left|2\eps\beta(x)\,\big(u_1 \bar u_2 + \bar u_1 u_2\big)\right| \le 2\eps|\beta(x)| \,\|U(x)\|^2\,.
$$
This implies \eqref{U-est}. The estimate \eqref{Uex-bound} is now a simple consequence of  \eqref{U-est}, presupposing that one proves some $\eps$-independent bounds on the initial condition $\|U_{ex}(x_d)\|$, or equivalently $\|(\chi_{ex}(x_d),\,1)^\top\|$. The latter norm is clearly bounded below by 1, and it is also bounded above due to the a-priori estimate on $\|\chi_{ex}\|_{C[0,x_d]}$ from Lemma \ref{lem3.1}. 
Hence there exist constants $0<C_1,\,C_2$, independent of $0<\eps\le\eps_0$, such that
\be\label{IC-bound}
  C_1\le \|U_{ex}(x_d)\| \le C_2\,,\qquad \forall 0<\eps\le\eps_0\,,
\ee
leading to \eqref{Uex-bound}.
\finproof

\subsection{Review of the WKB-marching method for the oscillatory region}\label{SEC331}
In this subsection we shall first review the WKB-marching method for solving the IVP \eqref{region2_ev} (or, equivalently, \eqref{EQU}). Then we recall its error estimates from \cite{ABAN11}.

Following  \cite{ABAN11} this method consists of two parts, first an analytic transformation of \eqref{region2_ev} or \eqref{EQU} into a less oscillatory problem, 
and second the discretization of the smooth problem on a coarse grid in an $\eps$-uniform manner. As shown in \cite{ABAN11}, the analytic WKB-transformation reviewed here is related to using oscillatory WKB-functions of second order, $\varphi_2^{os}(x)$. \\

\noindent
\underline{Part 1 -- analytic transformation:} The starting point is the vectorial IVP \eqref{EQU}. 
The vector function $U\in\CC^2$ is then transformed to the new unknown $Z\in\CC^2$ by
$$Z(x) =\left(
\begin{array}{c}
 z_1\\[2mm]
\ds z_2
\end{array}
\right) 
:= \exp \left(-{i\over \eps} \Phi^\eps(x) \right) P\,U(x)\,, \quad \forall x \in [x_d,1]\,,
$$ 
with the matrices
$$
P := {1\over \sqrt{2}}
\left(
\begin{array}{cc}
i&1\\
1&i
\end{array}
\right)\quad; \quad \Phi^\eps (x) := 
\left(
\begin{array}{cc}
\ds \phi^{\eps}(x) &0\\
0&-\phi^{\eps}(x) \end{array}
\right)\,,
$$
and the (real valued) phase function 
\begin{equation} \label{phase}
\phi^{\eps}(x):=\int_{x_d}^x \left( \sqrt{a(y)} - \eps^2 \beta(y)\right) \,dy\,.
\end{equation}
This change of unknown leads to the smooth ODE-system
\begin{equation} \label{EQZ}
\left\{
\begin{array}{l}
\ds {dZ\over dx} = \eps N^\eps Z\,,\quad x_d<x<1\,,\\[3mm]
\ds Z(x_d):=P\,U_{ex}(x_d)\,.
\end{array}
\right.
\end{equation}
Here, the $2\times2$--matrix function
$$
N^\eps(x) := \beta(x)\, \exp ( -{i\over \eps} \Phi^\eps )\, \left(
\begin{array}{cc}
0&1\\
1&0
\end{array}
\right) \, \exp ( {i\over
\eps} \Phi^\eps )\,,$$
is 
bounded independently of $\eps$. It is off-diagonal, with the entries 
$$
N^\eps_{1,2} (x)= {\beta}(x) 
e^{-\frac{2i}{ \eps} \phi^{\eps}(x)}\,,\quad
N^\eps_{2,1} (x)= {\beta}(x) e^{\frac{2i}{ \eps} \phi^{\eps}(x)}\,.
$$
This finishes the analytical transformation, and the goal of the second part is to provide an $\eps$-uniform discretization of \eqref{EQZ} that is second order w.r.t.\ the mesh size.\\

\noindent
\underline{Part 2 -- numerical discretization:}
First we partition the interval $[x_d,1]$ into $x_d=x_N < x_{N+1} <  \cdots < x_M = 1$. As in \S\ref{SEC32} we denote the mesh size by $h_n:=x_{n+1}-x_n$ as well as $h:= \max_{n=1, \cdots,M-1} \{h_n\}$.

With the initial condition $Z_N:=P\,U_N\in\CC^2$ and $U_N:=U_{ex}(x_d)\in\CC^2$ given, the marching scheme reads as follows (see \cite{ABAN11}):
\be\label{marching-scheme}
Z_{n+1} = (I + A_{n}^1+ A_{n}^2) \,Z_{n}\,,\qquad  n=N,...,M-1\,,
\ee
with the $2\times2$--matrices
$$
\begin{array}{lll}
\ds  A_{n}^1:=& \\[3mm]
&&\hspace{-20mm} - i \eps^2 \!\!\left(\!\!
\begin{array}{cc}
0&\hspace{-8mm} \beta_0(x_{n})e^{-{2i \over \eps} \phi(x_{n})} -\beta_0(x_{n+1})e^{-{2i \over \eps} \phi(x_{n+1})}\\[3mm]
\beta_0(x_{n+1}) e^{{2i \over \eps} \phi(x_{n+1})} -\beta_0(x_{n})e^{{2i \over \eps} \phi(x_{n})}
&\hspace{-8mm} 0
\end{array}
\!\!\!\right)\\[6mm]
&&\hspace{-20mm} \ds + \eps^3 \!\!\left(\!\!
\begin{array}{cc}
0&\hspace{-8mm} \beta_1(x_{n+1})e^{-{2i \over \eps} \phi(x_{n+1})} -\beta_1(x_{n})e^{-{2i \over \eps} \phi(x_{n})}\\[3mm]
\beta_1(x_{n+1}) e^{{2i \over \eps} \phi(x_{n+1})} -\beta_1(x_{n})e^{{2i \over \eps} \phi(x_{n})}
&\hspace{-8mm} 0
\end{array}
\!\!\right)\\[6mm]
&&\hspace{-20mm} \ds + i \eps^4 \beta_2(x_{n+1})  \left(
\begin{array}{cc}
0& -e^{-{2i \over \eps} \phi(x_{n})} H_1(-{2 \over \eps} S_{n})\\[3mm]
e^{{2i \over \eps} \phi(x_{n})} H_1({2 \over \eps} S_{n})&0
\end{array}
\! \right)\\[6mm]
&&\hspace{-20mm} \ds -  \eps^5 \beta_3(x_{n+1})  \left(
\begin{array}{cc}
0& e^{-{2i \over \eps} \phi(x_{n})} H_2(-{2 \over \eps} S_{n})\\[3mm]
e^{{2i \over \eps} \phi(x_{n})} H_2({2 \over \eps} S_{n})&0
\end{array}
\! \right)\,,
\end{array}
$$
$$
\begin{array}{lll}
\ds  A_{n}^2&:=& \ds - i \eps^3 (x_{n+1} -x_n) { \beta(x_{n+1}) \beta_0(x_{n+1}) +\beta(x_{n}) \beta_0(x_{n}) \over 2} 
\left(
\begin{array}{cc}
\ds 1&0\\
\ds  0&\ds -1
\end{array}
\right)\\[6mm]
&& \ds - \eps^4 \beta_0(x_{n}) \beta_0(x_{n+1}) 
\left(
\begin{array}{cc}
\ds H_1(-{2\over \eps}S_{n})&0\\[5mm]
\ds 0& \ds  H_1({2\over \eps}S_{n})
\end{array}
\right)\\[7mm]
&& \ds +i  \eps^5\beta_1(x_{n+1})  [\beta_0(x_{n})-\beta_0(x_{n+1})] \left( 
\begin{array}{cc}
\ds H_2(-{2\over \eps}S_{n})
&\ds 0\\[4mm]
\ds 0
&\ds -H_2({2\over \eps}S_{n})
\end{array}
\right)\,.
\end{array}
$$ 
Here we used the notation
$$
\beta_0(y) := \frac{\beta}{2 (\sqrt a -\eps^2\beta)} (y)\,; 
\qquad \beta_{k+1}(y) := \frac{1}{2 \phi'(y)} {d\beta_k \over dy}(y)\,,\quad k=0,\,1,\,2,
$$
$$
H_1(\eta)  := e^{i\eta} - 1\,, \qquad H_2(\eta)  := e^{i\eta} - 1-i\eta\,,
$$
and the discrete phase increments
$$
S_{n}\! := \phi(x_{n+1}) - \phi(x_{n})
=\! \int_{x_n}^{x_{n+1}} \!\!\! \left(\sqrt{a(y)}-\eps^2 \beta(y)\right) dy\,.
$$
Remark that for notational reasons we omitted in the aforementioned description of the scheme the $\eps$-index.
Furthermore we assumed that the two functions $\phi$ and $\beta$ (the latter involving the derivatives $a'$, $a''$) are explicitly ``available''. Alternatively, $\phi$, $a'$ and $a''$ could be approximated numerically. But, for simplicity, we shall not include such errors in the subsequent error analysis.

Finally we have to transform back to the $U$-solution vector via
\be \label{Transfo_ZU}
U_n=P^{-1} \,e^{{i \over \eps} \Phi^{\eps} (x_n)} \,Z_n\,, \qquad n=N+1,...,M\,,
\ee
which concludes the review of the WKB-marching algorithm.
\medskip

The following lemma is the discrete analogue of Lemma \ref{U-estimate}.
\vspace{0.2cm}
\begin{lemma}\label{lem:real-preservation}
Let Hypothesis B hold and let the initial condition $U_N\in\RR^2$. Then the iteration \eqref{marching-scheme}-\eqref{Transfo_ZU} determines a well-defined sequence satisfying $U_{n}\in\RR^2$ $\forall n=N+1,...,M$. Furthermore $\exists\, \tilde\eps_0\in(0,\eps_0]$ such that
\be\label{Un-stability}
  C_5\le \|U_n\| \le C_6\,, \qquad n=N,...,M\,,
\ee
with some positive constants $C_5,\,C_6$ that are independent of $0<\eps\le \tilde\eps_0$ and the numerical grid on $[x_d,1]$.
\end{lemma}\\
\vspace{0.2cm}

\debproof
Let us start by analysing the propagation matrix $B_n:=I+A_n^1+A_n^2\in\CC^{2\times2}$ of the vector $Z_n$ as defined in \eqref{marching-scheme}. A straightforward computation reveals its symmetry (which was also used in the proof of Proposition 3.3, \cite{ABAN11}):
$$
  \overline{(B_n)_{11}} = (B_n)_{22}\,,\quad \overline{(B_n)_{12}} = (B_n)_{21}\,,\qquad  n=N,...,M-1\,.
$$
This symmetry carries over to the matrix 
$$
  \tilde B_n =  \left( 
\begin{array}{cc}
b_1 &\ds b_2\\[4mm]
\bar b_2
&\bar b_1
\end{array}
\right) := e^{{i \over \eps} \Phi^{\eps} (x_{n+1})}\,B_n\,e^{-{i \over \eps} \Phi^{\eps} (x_n)} \,. 
$$
With this notation, the propagation matrix for the vector $U_n$ reads (cf. \eqref{marching-scheme}, \eqref{Transfo_ZU}):
$$
  P^{-1}\, \tilde B_n\, P = \left( 
\begin{array}{cc}
\Re b_1 + \Im b_2 & \Im b_1 + \Re b_2 \\[4mm]
-\Im b_1 + \Re b_2
& \Re b_1 - \Im b_2
\end{array}
\right)\in\RR^{2\times2}\,,
$$
where we used
$$
  P^{-1}= {1\over \sqrt{2}}
\left(
\begin{array}{cc}
-i&1\\
1&-i
\end{array}
\right)\,.
$$
This shows that $U_n\in \RR^2$.\\
Coming now to the bounds of $U_n$, a simple Taylor expansion for the matrices in \eqref{marching-scheme} yields $\|A_n^1\|\le C\eps h_n$,  $\|A_n^2\|\le C\eps^3 h_n$, and hence with some constant $C_7>0$:
$$
  \| Z_{n+1}\| \le \| Z_{n}\| (1+C_7 \eps h_n) \le \| Z_{n}\| e^{C_7 \eps h_n} \le \| Z_N \| e^{C_7 \eps (1-x_d)} 
  \le C_2 \| P\| e^{C_7 \eps (1-x_d)}\,,
$$
where we used $Z_N:=P\,U_{ex}(x_d)$ and the estimate \eqref{IC-bound} in the last step.
Next we consider the transformation matrices in \eqref{Transfo_ZU}: $P^{-1}$ is independent of $\eps$ and $h$, and $e^{{i \over \eps} \Phi^{\eps} (x_n)}$ is unitary. This implies the upper bound in \eqref{Un-stability}.

For the lower bound we choose $\tilde\eps_0:=\min\{\eps_0,\,\frac{1}{2C_7(1-x_d)}\}>0$ such that
$$
  \|A_n^1+A_n^2\|\le C_7\eps h_n\le 0.5\,;\qquad \forall 0<\eps\le \tilde\eps_0\,,\quad \forall 0<h_n\le 1-x_d\,.
$$
With the elementary estimate $\frac{1}{1-y}\le4^y$ for $y\in[0,0.5]$ we then obtain
$$
  \|(I+A_n^1+A_n^2)^{-1}\| \le \frac{1}{1-\|A_n^1+A_n^2\|} \le 4^{\|A_n^1+A_n^2\|}\,.
$$
This allows to estimate the backwards propagation $Z_n=(I+A_n^1+A_n^2)^{-1}\,Z_{n+1}$ as
$$
   \| Z_{n+1}\| \ge \| Z_{n}\|\, 4^{-C_7 \eps h_n} \ge \| Z_{N}\|\, 4^{-C_7 \eps (1-x_d)}\,,
$$
and the lower bound on $\|U_n\|$ follows as before.
\finproof\\
Due to the above lemma we have to restrict the range of admissible $\eps$--values:\\

{\bf Hypothesis B'} {\it Let the assumptions of Hypothesis B hold, but with $\eps_0$ replaced by $\tilde\eps_0$ from Lemma \ref{lem:real-preservation}.}
\subsubsection{Error and stability estimates for the WKB-marching method}

In this subsection we recall the main Theorem 3.1 of \cite{ABAN11}, providing error bounds for the marching method \eqref{marching-scheme}-\eqref{Transfo_ZU}, used for solving the IVP \eqref{region2_ev} or, equivalently, \eqref{EQU}.\\

\begin{theorem} \label{THM_princ} {\bf (Convergence WKB-IVP)}
Let Hypothesis B be satisfied and let $U_{ex}(x)$ denote the exact solution to the IVP \eqref{EQU}. Then the global error of the 
second order scheme \eqref{marching-scheme}-\eqref{Transfo_ZU} satisfies
\be \label{error_Z_2ORD}
\qquad  \|U_{ex}(x_n)-U_n\|_{} \le C {h^{\gamma}
\over \eps} +C \eps^3 h^2\,, \; N\le n\le M\,,
\ee 
with $C$ independent of $n$, $h$, and $\eps$.
Here, $\gamma >0$ is the order of the chosen numerical integration
method for computing the phase integral 
\be \label{PH}
\Phi^{\eps}(x)= \int_{x_d}^x \left(\sqrt{a(y)} -\eps^2 \beta(y) \right) dy
\, \left(
\begin{array}{cc}
1&0\\0&-1
\end{array}
\right)
\,
\ee
for the back-transformation \eqref{Transfo_ZU}.
\end{theorem}\\
\vspace{0.2cm}

We remark that the term ${h^{\gamma} \over \eps}$ of \eqref{error_Z_2ORD} may or may not be present in real computations, depending on the chosen coefficient function $a(x)$. If $a(x)$ is piecewise linear or piecewise quadratic, e.g., the phase integral $\phi^\eps(x)$ can be computed analytically. Hence, this term would not appear in such cases. In the numerical tests performed in \S\ref{SEC4} below, we shall only consider such examples of exactly computable phase functions and shall hence not include this error term in the error analysis of \S\ref{SEC34}.
\subsection{Convergence results for the overall hybrid WKB method} \label{SEC34}

In this section we shall combine the error analysis of the previous two sections and adapt it to the algorithm for coupling two regions. 
To this end we have to fix the numerical analogues of the continuous coupling conditions in \eqref{region2_ev}. First we shall (of course) use $\varphi_h(x_d):=\chi_h(x_d)$. But for the initial condition of the derivative there are two options, namely $\eps\varphi_h'(x_d):=\eps\chi_h'(x_d)$ or $\eps\varphi_h'(x_d):=1$ (taken from the exact value in \eqref{region1_ev}). We shall use the second option for the following reasons: On the one hand it avoids the numerical error of $\chi_h'(x_d)$, where we recall that $\chi_h'$ is discontinuous (and hence less accurate) at the grid points.  And on the other hand, this choice will facilitate the a-priori estimate needed for the scaling in Step 3.

Since the numerically used initial data $\chi_h(x_d)$ deviates from its exact value $\chi_{ex}(x_d)$, this gives rise to an additional error component to be considered: Let thus $\hat U_{ex}(x)$ denote the exact solution to the ODE \eqref{EQU}, but with the following perturbed 
initial condition:
\bea\label{numIC}
  \hat U_{ex}(x_d):=
  A(x_d+)\left( \begin{array}{c} 
  \chi_{h}(x_d) \\ 1
  \end{array} \right)\,.
\eea
Using the a-priori estimate \eqref{U-est} leads to the error
\bea\label{error1}
   && \|U_{ex}-\hat U_{ex}\|_{C[x_d,1]} \le
   \left\|A(x_d+)\,\big(e_h(x_d),\,0\big)^\top\right\|\,\exp\Big[ \eps \int_{x_d}^1  |\beta(y)| dy\Big] \nonumber\\
   &&\qquad\qquad\quad = \Big[ a(x_d+)^{1/2}+\frac{\eps^2}{16}a^{-5/2}(x_d+) a'(x_d+)^2\Big]^{1/2}
   \,|e_h(x_d)|
   \exp\Big[ \eps \int_{x_d}^1  |\beta(y)| dy\Big]\,, 
\eea
where we used the notation $e_h=\chi_{ex}-\chi_h$.\\

The following convergence analysis of the hybrid method uses several different solution functions (exact, numerical, etc.). To keep the notation straight we summarize it in the following table, both for the evanescent and oscillatory regions. The superscript ${}{}^{(}{}'{}^{)}$ signifies that we refer to both the function and its first derivative. 

\[
\begin{array}{|r|l|}
\hline 
\hbox{\rm evanescent} &   \\ \hline
\phantom{I^{I^{I^I}}} 
\chi_{ex}\!\!\!{}^{(}{}'{}^{)}(x) & \hbox{\rm exact solution of Step 1, eq. \eqref{region1_ev}, \eqref{VF_ev}} \\
\chi_{h}\!{}^{(}{}'{}^{)}(x) & \hbox{\rm numerical solution of Step 1, eq. \eqref{VF_ev_h}} \\ \hline
\hbox{\rm oscillatory} &   \\ \hline
\phantom{I^{I^{I^I}}} 
\varphi_{ex}\!\!\!{}^{(}{}'{}^{)}(x) & \hbox{\rm exact solution of Step 2, with exact IC $\chi_{ex}(x_d)$, eq. \eqref{region2_ev}} \\
U_{ex}(x) & \hbox{\rm exact solution vector of Step 2, with exact IC $\chi_{ex}(x_d)$, eq. \eqref{EQU}} \\
\hat U_{ex}(x) & \hbox{\rm exact solution vector of Step 2, with numerical IC $\chi_{h}(x_d)$, eq. \eqref{EQU}, \eqref{numIC}} \\
U_n & \hbox{\rm numerical solution vector of Step 2, with numerical IC $\chi_{h}(x_d)$, eq. \eqref{marching-scheme}-\eqref{Transfo_ZU}} \\
\hline
\hbox{\rm hybrid} &   \\ \hline
\phantom{I^{I^{I^I}}} 
\psi_{ex}\!\!\!{}^{(}{}'{}^{)}(x) 
& \hbox{\rm exact solution after scaling in Step 3, eq. \eqref{SchBVP}, \eqref{psi-scaling_ev}} \\
\psi_{h}\!{}^{(}{}'{}^{)}(x)\,;\; \psi_{h,n}  & \hbox{\rm numerical solution after scaling in Step 3 (on $[0,x_d]\, ;\; [x_d,1]$)} \\
\tilde U_{ex}(x) \mbox{ on } [x_d,1] 
& \hbox{\rm exact solution vector after scaling in Step 3, with exact IC} \\
\tilde U_n  \mbox{ on } [x_d,1] 
& \hbox{\rm numerical solution vector after scaling in Step 3, with numerical IC} \\
\hline
\end{array}
\]

For clarity, we summarize here also the numerical analogue of the three steps in \eqref{region1_ev}-\eqref{alpha-eq_ev}, referring to the two regions in Fig.\ \ref{2regions}: \\

\noindent
\underline{Step 1 -- WKB-FEM for $\chi_h$ in region (1):} Find $\chi_h\in\mathcal V_h$ solving
\be \label{FEM-region1}
  b(\chi_h,\theta_h) = L(\theta_h)\,, \quad \forall \theta_h\in\mathcal V_h\,.
\ee
This yields $\chi_h(x_d)\in\RR$ with an error $|e_h(x_d)| \le C\sqrt h\min\{\eps^{3/2},h^{3/2}\}$ (see Thm. \ref{EV_conv}).
\medskip

\noindent
\underline{Step 2 -- WKB-marching method for $\varphi_h$ in region (2):}
As initial condition for the \linebreak
marching scheme we use $U_N:=\hat U_{ex}(x_d)\in\RR^2$ given by \eqref{Un-stability}. Applying the scheme \eqref{marching-scheme}--\eqref{Transfo_ZU} iteratively we compute the vectors $U_n=(u_n^1,\,u_n^2)^\top\in\RR^2; \,n=N+1, \cdots,M$.
\medskip

\noindent
\underline{Step 3 -- Scaling of the auxiliary wave functions $\chi_h,\,\varphi_h$:}
\be \label{psi-scaling_num}
\qquad \psi_h(x):=\left\{
\begin{array}{l}
\ds \tilde\alpha\,\chi_h(x)\,, \qquad\qquad\qquad\qquad x \in [0,x_d)\,, \\[3mm]
\ds \tilde\alpha\,\varphi_{h,n}=\tilde\alpha\,u_n^1\,a(x_n)^{-1/4}\,, \quad x\in\{x_N,...,x_M\} \,,
\end{array}
\right.
\ee
$$
\psi_h'(x):=\left\{
\begin{array}{l}
\ds \tilde\alpha\,\chi_h'(x)\,, \qquad\qquad\qquad\qquad x \in [0,x_d)\,, \\[3mm]
\ds \tilde\alpha\,\varphi_{h,n}'=\tilde\alpha \Big[\frac{u_n^2}{\eps}\,a(x_n)^{1/4}
-\frac{a'(x_n)}{4}\, a(x_n)^{-5/4}\,u_n^1\Big]\,, \quad x\in\{x_N,...,x_M\} \,,
\end{array}
\right.
$$
with the scaling parameter $\tilde\alpha\in\CC$ defined in analogy to \eqref{alpha-eq_ev}:
\be\label{alpha-scaling}
  \tilde\alpha = \tilde\alpha (u_M^1,\,u_M^2) := \frac{-2\i a(1)^{1/4}}{u_M^2-\big[i+\frac\eps4 a(1)^{-3/2}a'(1)\big] u_M^1} \,.
\ee
The statement \eqref{psi-scaling_num} reveals that our final numerical solution $\psi_h$ is continuous in the evanescent region, but discrete in the oscillatory region.
Note also that the relation between $U=(u_1,\,u_2)^\top$ and $(\varphi,\,\eps \varphi')^\top$, given by \eqref{TRANSFO}, provides a connection between the two scaling functions $\alpha$ and $\tilde\alpha$, {\it i.e.}
$$
  \tilde\alpha(u_1,\,u_2) = \alpha(\varphi,\,\varphi')\,.
$$
Let us finally also recall that the solution $\chi_h$ as well as the vector $U_n\ne0$ of Step 2 are real valued. The final (numerical) solution $\psi_h$ only becomes complex valued due to the multiplication by $\tilde\alpha$ in Step 3. Note that the denominator of \eqref{alpha-scaling} cannot vanish for $U\in \RR^2\setminus \{0\}$, which makes the scaling well defined. This map $\tilde\alpha$ satisfies moreover the following simple properties:
\vspace{0.2cm}

\begin{lemma}\label{lem:Lipschitz}
For each fixed $\delta>0$, the ($\eps$-dependent) map $\tilde \alpha:\,U\in\RR^2\setminus B_{\delta}(0) \to \CC$ is Lipschitz continuous with some constant $L>0$ and bounded by some constant $C_8$. Both constants can be chosen uniformly w.r.t.\ $0<\eps\le\eps_0$ and are dependent on $\delta$.
\end{lemma}
\medskip

The following error analysis of the hybrid scheme is the main result of this paper.\\

\begin{theorem}\label{th:hybrid-conv} {\bf (Convergence WKB-hybrid)}
Let Hypotheses A and B' be satisfied. Then $\psi_h$, the numerical solution to the hybrid scheme \eqref{FEM-region1}-\eqref{alpha-scaling}, satisfies the following error estimates, compared to the exact solution $\psi_{ex}$ of the algorithm \eqref{region1_ev}-\eqref{alpha-eq_ev}:
\begin{enumerate}
\item[a)] In the evanescent region $[0,x_d)$ we have 
\bea\label{hybrid-evan}
  \|\tilde e_h\|_{L^2(0,x_d)} &\le& C \sqrt{\eps}\, h\, \min \{ \eps,h\}\,, \qquad\quad 
  \eps\,\|\tilde   
  e_h'\|_{L^2(0,x_d)} \le C \eps^{3/2}\, h\,,\nonumber\\
  \|\tilde e_h\|_{C[0,x_d]} &\le& C \,\sqrt{h}\, \min \{ \eps^{3/2},h^{3/2}\}\,,\quad 
  \eps  \|\tilde
  e_h'\|_{L^{\infty}(0,x_d)}\le C \eps\, \sqrt{h}\,\min \{ \sqrt{\eps},\sqrt{h}\} \,,\nonumber\\
  &&
\eea
with the notation $\tilde e_h(x):=\psi_{ex}(x)-\psi_h(x)$.
\item[b)] In the oscillatory region $[x_d,1]$ we have 
\be\label{hybrid-oscil}
  |\tilde e_{h,n}| +\eps |\tilde e_{h,n}'|  \le C \sqrt{h}\, \min \{ \eps^{3/2},h^{3/2}\}\,;
  \qquad n=N,...,M\,,
\ee
with the notation $\tilde e_{h,n}:=\psi_{ex}(x_n)-\psi_{h,n}$ and $\tilde e_{h,n}':=\psi_{ex}'(x_n)-\psi_{h,n}'$.
\item[c)] For the overall hybrid method one has then, over $[0,1]$, the estimates
$$
  \|\tilde e_h\|_\infty \le C \,\sqrt{h}\, \min \{ \eps^{3/2},h^{3/2}\}\,,\qquad \eps  \|\tilde
  e_h'\|_\infty \le C \eps\, \sqrt{h}\,\min \{ \sqrt{\eps},\sqrt{h}\} \,,
$$
with the notation $\|\tilde e_h\|_\infty := \max\{ \|\tilde
  e_h\|_{L^{\infty}(0,x_d)}\,;\; \ds \max_{n=N,...,M} |\tilde e_{h,n}| \}$.
\end{enumerate}
\end{theorem}
\medskip

Note that the error of $\psi_h$ is globally of the order $\mathcal O(h^2)$. The derivative $\psi_h'$ is also correct to order $\mathcal O(h^2)$ in the oscillatory region, but only $\mathcal O(h)$ in the evanescent region. This can be explained as follows: The $\mathcal O(h)$--error of $\chi_h'$ is not propagated to the oscillatory region, as we use the exact value of the derivative in the matching condition \eqref{numIC}. Moreover, the scaling step does not change the error orders in the evanescent region.\\

\debproof[of Theorem \ref{th:hybrid-conv}]
Statement a) is a consequence of Thm. \ref{EV_conv}. Let us then  continue by estimating the error of the numerical solution $U_n$ compared to the exact solution $U_{ex}(x)$, i.e. prior to the scaling Step 3: Using \eqref{error1}, \eqref{error_Z_2ORD}  and Theorem \ref{EV_conv} we obtain
\bea\label{error2}
  \|U_{ex}(x_n)-U_n\| &\le& \|U_{ex}(x_n)-\hat U_{ex}(x_n)\| + \|\hat U_{ex}(x_n)-U_n\| \nonumber\\
  &\le& C\,|e_h(x_d)| + C\eps^3h^2\\
  &\le& C\sqrt h \min\{\eps^{3/2},\,h^{3/2}\} 
  \,;\quad n=N,...,M\,.\nonumber
\eea
We continue with estimating the error propagation due to the non-linear scaling in Step 3. Due to Lemma \ref{lem:Lipschitz}, the map $\tilde \alpha:\,U\in\RR^2\setminus B_{\min\{C_3,C_5\}}(0) \to \CC$ is Lipschitz continuous with some constant $L>0$ and bounded by some constant $C_8$. Both of these constants can be chosen independent of $0<\eps\le\eps_0$, as the choice $\delta:=\min\{C_3,C_5\}$ for the domain of $\tilde\alpha$ uses the lower bounds on $U_{ex}$ from \eqref{Uex-bound} and on $U_n$ from Lemma \ref{lem:real-preservation}. Here it is crucial that both the exact solution vector $U_{ex}(1)$ and the numerical solution vector $U_M$ have real components. Then, Lemma \ref{lem:Lipschitz} shows that $\tilde\alpha$ is Lipschitz and bounded on $\RR^2\setminus B_\delta(0)$.\\

\noindent
\underline{Part a)} For the evanescent region $[0,x_d)$ we estimate the difference between the exact solution and the numerical solution (both after scaling)
$$
  \psi_{ex}(x) =\tilde\alpha(U_{ex}(1))\,\chi_{ex}(x)\,, \quad
  \psi_{h}(x) = \tilde\alpha(U_M)\,\chi_{h}(x)\,.
$$
This yields
\bean
  |\psi_{ex}(x) - \psi_{h}(x)|
  &\le& |\tilde\alpha(U_{ex}(1))-\tilde\alpha(U_M)|\;|\chi_{ex}(x)|
  +|\tilde\alpha(U_M)|\;|\chi_{ex}(x)-\chi_h(x)| \\
  &\le& L\;\|U_{ex}(1)- U_M\|\;|\chi_{ex}(x)| + C_8\; |e_h(x)|\\
  &\le& C\sqrt h \min\{\eps^{3/2},\,h^{3/2}\}\;|\chi_{ex}(x)| + C_8\; |e_h(x)|\,,
\eean
where we used \eqref{error2} in the last step. Using now Lemma \ref{lem3.1} to estimate $\chi_{ex}\!\!\!{}^{(}{}'{}^{)}(x)$ and Theorem \ref{EV_conv} for $e_h\!{}^{(}{}'{}^{)}(x)$ yields the four 
error estimates of \eqref{hybrid-evan}.\\

\noindent
\underline{Part b)} For the oscillatory region $[x_d,1]$ we estimate the difference between the exact solution and the numerical solution (both after scaling)
$$
  \psi_{ex}(x_n) =\tilde\alpha(U_{ex}(1))\,\varphi_{ex}(x_n)\,, \quad 
  \psi_{h,n} = \tilde\alpha(U_M)\,\varphi_{h,n}\,,
$$
at the grid points $x_n;\,n=N,...,M$. Here it is again more convenient to use the vector notation from \eqref{def-U}, where we introduce the notations $\tilde U_{ex}(x) := \tilde \alpha U_{ex}(x)$ and $\tilde U_{n} := \tilde \alpha U_{n}$ for the exact and, respectively, numerical solution after scaling in Step 3:
\bean
  \|\tilde U_{ex}(x_n) - \tilde U_{n}\| 
  &=& \|\tilde\alpha(U_{ex}(1)) \,U_{ex}(x_n) - \tilde\alpha(U_M)\,U_{n}\| \\
  &\le& |\tilde\alpha(U_{ex}(1))-\tilde\alpha(U_M)|\;\|U_{ex}(x_n)\|
  +|\tilde\alpha(U_M)|\;\|U_{ex}(x_n)- U_n\| \\
  &\le& L\;\|U_{ex}(1)- U_M\|\;C_4 + C_8\; \|U_{ex}(x_n)- U_n\|\\
  &\le& C\sqrt h \min\{\eps^{3/2},\,h^{3/2}\} \,,
\eean
where we used \eqref{Uex-bound} in the penultimate line, and \eqref{error2} twice in the last line. 
Using the norm equivalence then yields the error estimate \eqref{hybrid-oscil}.\\

\noindent
\underline{Part c)} is just a combination of the previous two parts.
\finproof

\section{Numerical tests of the hybrid WKB method} \label{SEC4}
\begin{figure}[ht!]
\begin{center}
 \includegraphics[scale=.43]{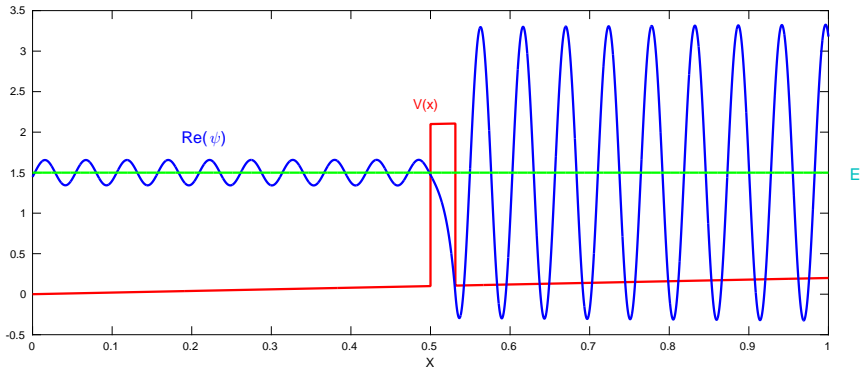}
 \caption{Tunnelling structure with injection of a plane wave from the right boundary $x=1$.
 Red curve: piecewise linear potential $V(x)$ with applied bias (since $V(1)>V(0)$). It is discontinuous at $x_c=0.5$ and $x_d=0.53125$.
 Blue curve: $\Re(\psi(x))$, real part of the wave function that is partly transmitted, but mainly reflected in this example; $\eps=0.01$.
 Green line: energy of the injected plane wave with $E<\max V(x)$. (colours only online)}
 \label{fig:4.1}
\end{center}
\end{figure}

The aim of this section is to present numerical results obtained with the WKB-coupling scheme introduced in Section \ref{SEC2} and to compare 
these results with the error analysis established in Section \ref{SEC34}. In particular, we present the results for 3 zones (oscillating-evanescent-oscillating, cf. \S\ref{SEC22}) corresponding to the passage or flow of electrons through a tunnelling structure (see Fig. \ref{fig:4.1}), with a piecewise linear and, respectively, piecewise quadratic potential $V(x)$, chosen such that the phase $\phi^\eps(x)$ is explicitly calculable. The reason for such a choice is to avoid having to care about the ${h^\gamma \over \eps}$ --error term in \eqref{error_Z_2ORD}, yielding hence an asymptotically correct scheme for fixed $h>0$ and $\eps \rightarrow 0$. 

\medskip
\noindent
\underline{Example 1:}
We start with the piecewise linear potential graphed in Fig. \ref{fig:4.1}. Note the small applied bias with $V(0)=0$, $V(1)=0.2$.

\begin{figure}[htbp]
\begin{center}
\includegraphics[width=6cm]{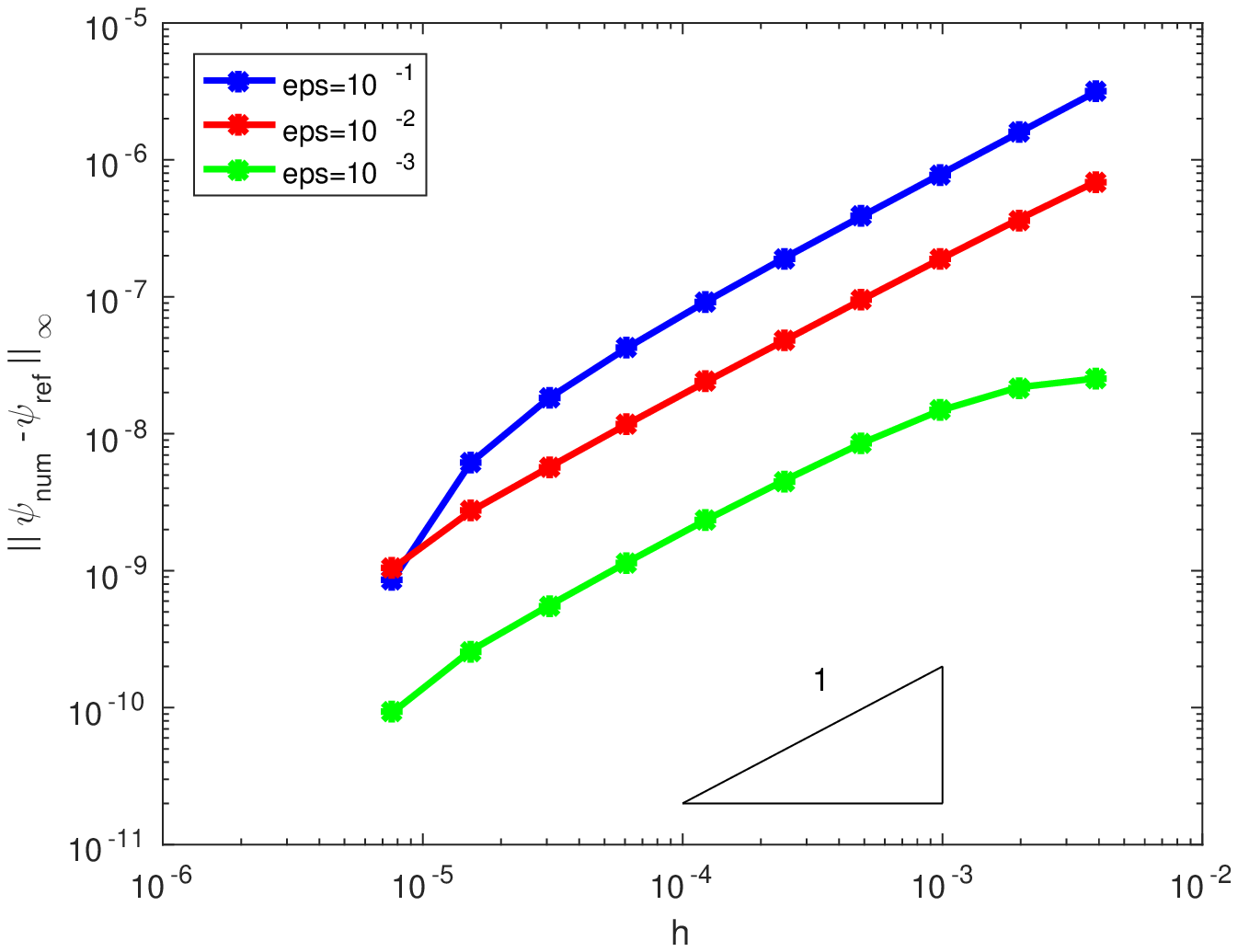}\hfill
\includegraphics[width=6cm]{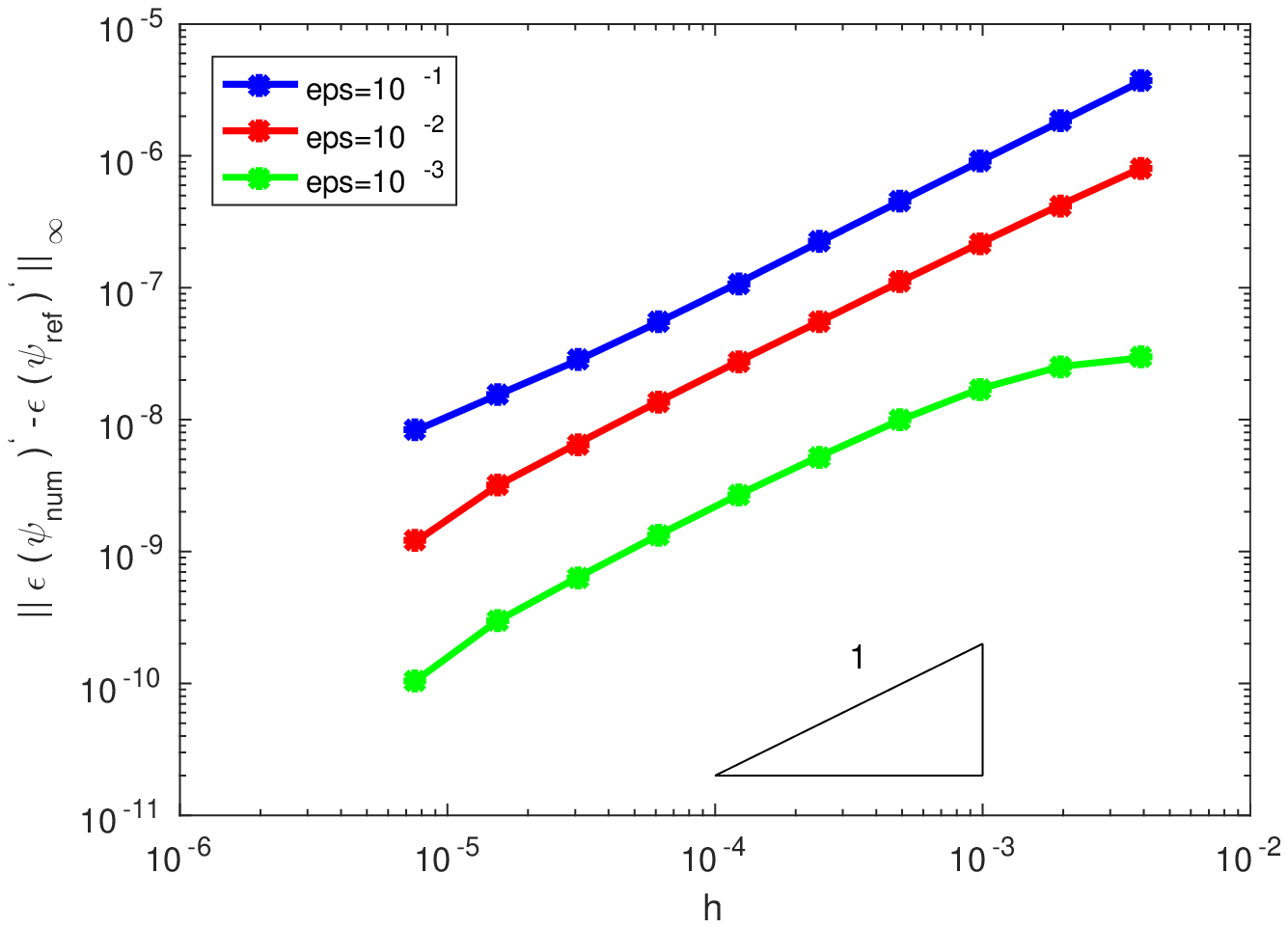}
\end{center}
\caption{\label{err_inf} {\footnotesize Absolute error (in the $L^\infty$-norm and $\log-\log$ scale) between the computed solution and a reference solution (obtained with $h=2^{-18}\approx 4\cdot 10^{-6}$), for the  piecewise linear potential from Fig. \ref{fig:4.1}. Left: $||\psi_{ref}-\psi_{num}||_\infty$. Right: $||\eps \psi'_{ref}-\eps \psi'_{num}||_\infty$}.}
\end{figure}

In Fig. \ref{err_inf} we plotted the numerical errors of the coupling method associated to the wave function $\psi$ (left figure) and to its derivative $\eps \psi'$ (right figure), as functions of the mesh size $h$ (in $\log-\log$ scale) and for three different $\eps$-values. In the oscillating regions we chose the second order method \eqref{marching-scheme}-\eqref{Transfo_ZU} and in the evanescent region the FEM \eqref{VF_ev_h}. The plotted errors are the $L^\infty$-errors between the numerical solution on the whole interval $[0,1]$ and a reference solution, computed with the same scheme but on a finer grid of $2^{18}$ points. It can be observed that the slopes in these two plots are approximatively one (for $h\gtrsim 3\cdot 10^{-5}$) and improving to $1.5$ for smaller values of $h$. For $\eps=0.1$ the slope of the $\psi$--error even improves up to 2 for the smallest values of $h$. This behaviour is in accordance with our numerical analysis in Theorem \ref{th:hybrid-conv}(c!
 ).

The $\eps$-dependence seems to be like $\mathcal O (\sqrt\eps)$ (for large values of $\eps$), improving to $\mathcal O (\eps)$ (for small values of $\eps$), and even $\mathcal O (\eps^{3/2})$ (for small values of $\eps$ and large $h$). Summarizing, the error of $\psi$ shows to be of order $\mathcal{O} (\min\{h^2,\,\sqrt\eps h^{3/2},\,\eps h,\,\eps^{3/2}\sqrt h \})$, which corresponds exactly to the estimates given in Theorem \ref{th:hybrid-conv}(c). The error of $\eps\psi'$ shows to be of order $\mathcal{O} (\min\{\sqrt\eps h^{3/2},\,\eps h,\,\eps^{3/2}\sqrt h \})$, which is even slightly better than the estimate from Theorem \ref{th:hybrid-conv}(c) (in the sense of including also an $\mathcal O (\sqrt\eps h^{3/2})$--behaviour). 

We mention that the obtained numerical errors are mainly those introduced by the WKB-FEM of the evanescent region. Indeed, in this evanescent region, the numerical error of the WKB-FEM is larger than the one obtained from the second order WKB marching method of the oscillating region (compare the estimates in the Theorems \ref{EV_conv} and \ref{THM_princ}). \\


\medskip
\noindent
\underline{Example 2:}
Next we consider a piecewise quadratic potential given by
$$
a(x):=c_1 (x + c_2)^2\,, \quad \forall x \in [0,x_c] \cup [x_d,1]\,; \quad a(x):=-c_1 (x + c_2)^2\,,\quad \forall x \in (x_c,x_d)\,,
$$
with $x_c=0.5$, $x_d=0.5+2^{-5}=0.53125$ and
$$
E=1.5\, , \quad V_1=V(1)=0.2\, , \quad c_2= -{ E + \sqrt{E^2 - V_1\, E} \over V_1}\,, \quad c_1= { E \over c_2^2}\,.
$$
\medskip 

\begin{figure}[htbp]
\begin{center}
\includegraphics[width=6cm]{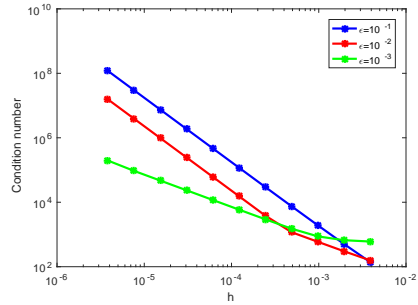}
\end{center}
\caption{\label{cond-number} {\footnotesize Condition number for the discrete BVP in Example 2, as a function of $h$, for three values of $\eps$.}}
\end{figure}

Before turning to the error plots we consider the condition number associated to solving the discrete variational problem \eqref{VF_ev_h} in the (intermediate) evanescent region. In Fig.\ \ref{cond-number} we plot this condition number as a function of $h$, for three different values of $\eps$. For $\eps=10^{-1},\,10^{-2}$ it grows like $\mathcal{O}(h^{-2})$ when $h\to0$, and for $\eps=10^{-3}$ it grows like $\mathcal{O}(h^{-1})$ (on the shown interval of $h$-values). We remark that this behaviour is not a problem in practice: For large $\eps$, the solution $\psi_{ex}$ is \emph{not} highly oscillatory and hence does not need a high spatial resolution. For small values of $\eps$, even a fine resolution would only lead to moderate condition numbers. 
Indeed, one observes a decrease of the condition number when $\eps$ gets smaller. This important feature is somehow related to the {\it asymptotic-preserving} property of the scheme.

Large condition numbers signify that the errors of numerical experiments also include significant contributions stemming from round-of errors and their accumulation. While the method-error (as estimated in Theorem \ref{th:hybrid-conv}(c)) decreases with decreasing $h$, the round-of errors could then increase in some situations, due to the increasing condition number. These arguments may lead to the idea that one cannot trust too much the reference solution in Fig. \ref{err_inf}, computed with $2^{18}$ points. In order to verify this suspicion, we decided to plot in the case of a piecewise quadratic potential two types of error curves.  

\begin{figure}[htbp]
\begin{center}
\includegraphics[width=6cm]{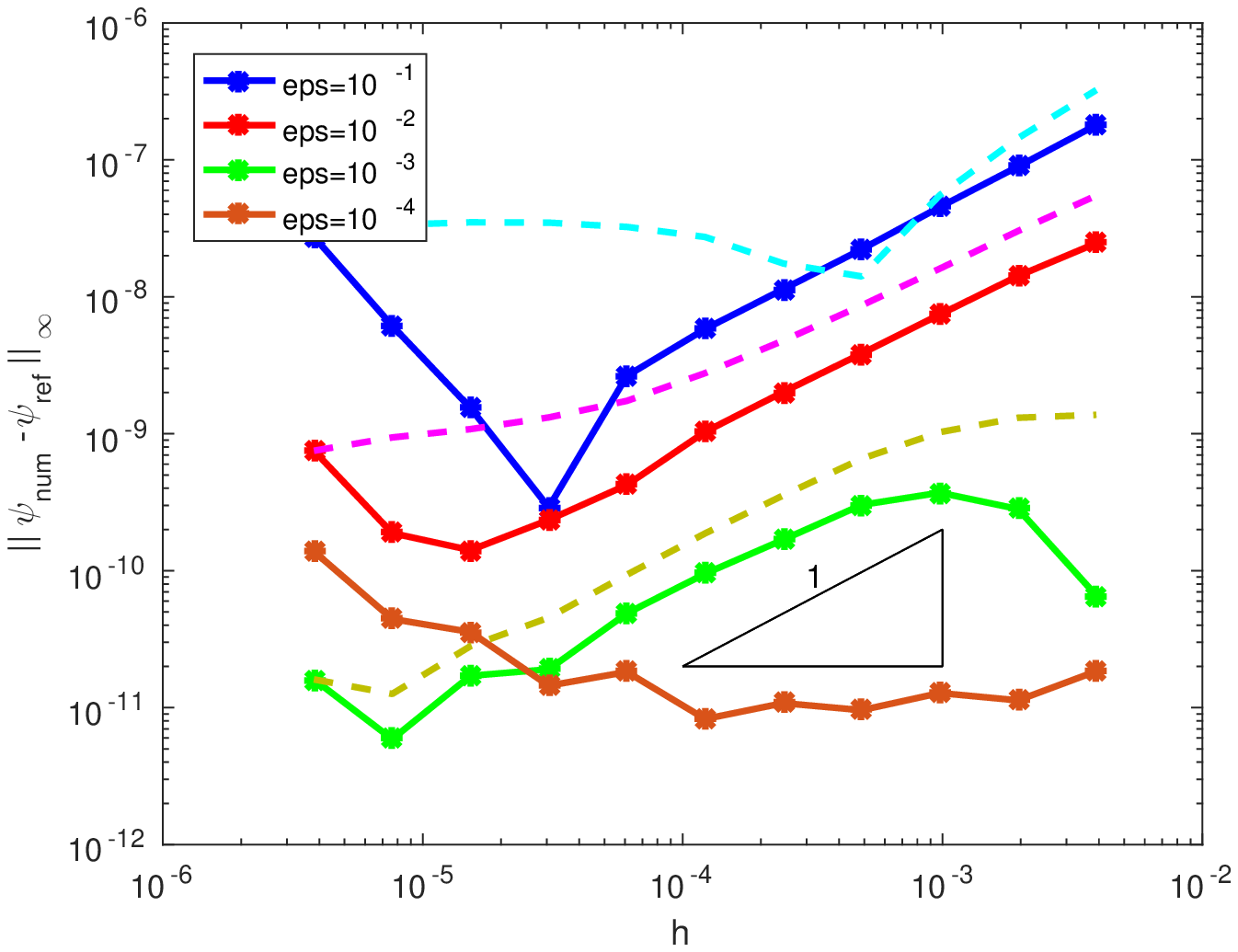}\hfill
\includegraphics[width=6cm]{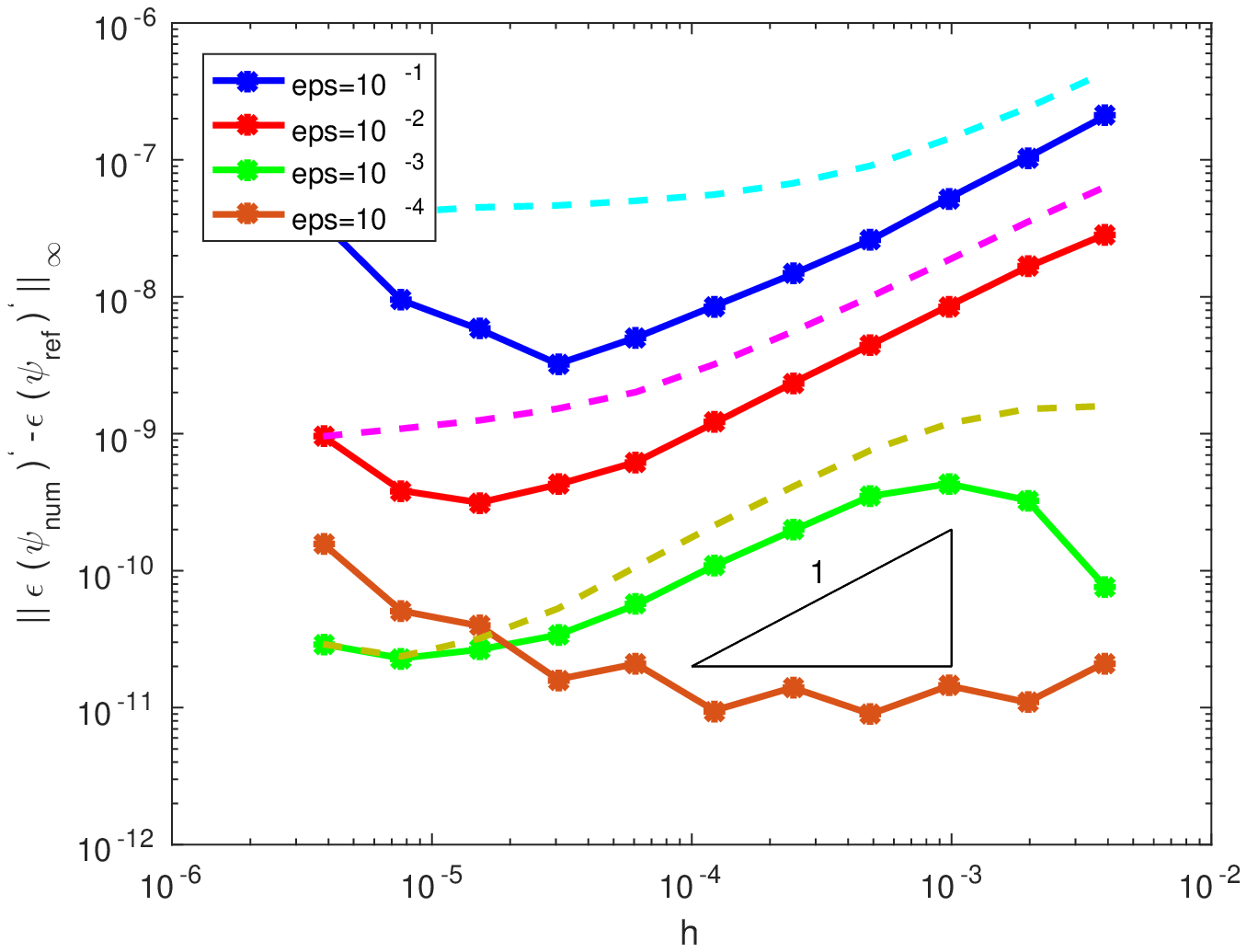}
\end{center}
\caption{\label{VQ_err_inf} {\footnotesize Absolute error (in the $L^\infty$-norm and $\log-\log$ scale) between the computed solution and a reference solution, for a piecewise quadratic potential. Left: $||\psi_{ref}-\psi_{num}||_\infty$. Right: $||\eps \psi'_{ref}-\eps \psi'_{num}||_\infty$}. For the dashed curves, $\psi_{ref}$ is computed with $h=2^{-19}$. The solid curves show the \emph{incremental error} when refining the mesh by a factor two. Left: $||\psi_{h_j}-\psi_{h_{j-1}}||_\infty$. Right: $||\eps \psi'_{h_j}-\eps \psi'_{h_{j-1}}||_\infty$.}
\end{figure}

To be more precise, in Fig.\ \ref{VQ_err_inf} we show the numerical errors of the wave function $\psi$ (left figure) and its derivative $\eps \psi'$ (right figure), as functions of the mesh size $h$ (in $\log-\log$ scale) and for four different $\eps$-values. The dashed lines correspond (as in Example 1) to the $L^\infty$-error between the numerical solution on the whole interval $[0,1]$ and a reference solution, computed with the same scheme but on a finer grid (here $h=2^{-19}$) whereas the solid lines correspond to the \emph{incremental error} when iteratively doubling the grid size, i.e. $||\psi_{h_j}-\psi_{h_{j-1}}||_\infty$ with $h_j=2h_{j-1}$. For a first order method, the former error is about twice as large as the latter (incremental) error. This can be understood from the geometric series of the incremental errors, since the  summands then differ by a factor of about 2. In Fig.\ \ref{VQ_err_inf} this difference is clearly visible for the solid red curves, pertainin!
 g to $\eps=10^{-2}$, and the 
corresponding dashed error curves (for large $h$).
The minimum of the incremental error (as a function of $h$) indicates the onset of significant round-of errors when reducing $h$. In Fig.\ \ref{VQ_err_inf} this is best visible for the solid blue and red curves, pertaining to $\eps=10^{-1},\,10^{-2}$. 

Furthermore remark that for $\eps=10^{-1},\,10^{-2},\,10^{-3}$ and $h\gtrsim 3\cdot 10^{-5}$, the error slopes are approximately one -- just like in Example 1. For smaller values of $h$ the error then gets polluted by round-of errors. For $\eps=10^{-4}$ the shown errors seem to be mostly due to round-of errors. They again increase for $h\lesssim 3\cdot 10^{-5}$.


\section{Turning points} \label{SEC5}
A \emph{turning point} of the Schr\"odinger equation \eqref{EQ_ref} is defined as a zero of the given coefficient function $a(x)$. Accordingly one also speaks about the \emph{order of a turning point}. We first remark that both error analyses, in \S\ref{SEC32} for the WKB-FEM and in \cite{ABAN11} for the WKB-marching method are \emph{not} valid for turning points. Therefore, we assumed in Hypothesis A and B that $a(x)$ is bounded away from zero. Furthermore, in the convergence Theorems \ref{EV_conv} and \ref{THM_princ} we did not keep track how the leading constant $C$ grows with $\tau_{ev},\,\tau_{os}\to0$. However, the failure of both WKB methods when approaching turning points also appeared in our numerical experiments (not included here).
The paper \cite{lub} considers a matrix generalization of our equation \eqref{EQ_ref}, but only for the oscillatory case. Their coefficient matrix $A(x)$ (generalizing our $a(x)$) is there assumed to be symmetric positive definite, satisfying the uniform lower bound $A(x) \ge \delta^2>0$. The proof of their Theorem 6.1 shows that their $L^\infty$--error bounds would grow like $\mathcal O(\delta^{-2}$).

This numerical failure can be understood easily: The WKB-ansatz \eqref{WKB}-\eqref{evWKB-fct} is not valid at turning points. In fact, close to a turning point of first order, solutions to \eqref{EQ_ref} are neither exponential nor oscillatory, but in a transition layer of thickness $\mathcal O(\eps^{2/3})$ they behave rather like Airy functions: Indeed, for $a(x):=-x$, a solution basis for \eqref{EQ_ref} is given by Ai$(\eps^{-2/3}x)$, Bi$(\eps^{-2/3}x)$, where Ai and Bi denote the Airy functions of first and second kinds. 

At a turning point the solution to \eqref{EQ_ref} clearly satisfies $\psi''=0$. This motivated to use linear FEM-ansatz functions in the numerical cell containing a turning point (cf. \S3.2.2 in \cite{Ne05}).

The quest for an appropriate replacement of the WKB-ansatz close to turning points has a long history in asymptotic analysis, starting with Langer \cite{La31}: For general coefficients $a(x)$ with a zero,
he found an asymptotic approximation for the solution of  \eqref{EQ_ref} that is valid uniformly in $x$, including the turning point. For a first order turning point his approximation is a composite function involving Airy functions and the phase function (like $\sigma(x)$ defined in \eqref{NOTZ}).
For details on first order turning points we refer to \S4.3 in \cite{Ho95}, and to \S7.3 in \cite{Ny73} for higher order turning points.

The above mentioned approximation formulas of Langer have also been used for numerical computations, mostly for Schr\"odinger eigenvalue problems \cite{GGG91,GGG92,SBW08}. In the physics literature, this strategy is frequently called \emph{Modified Airy function} (MAF) technique. It relies on evaluating the explicit formulas of approximate solutions, but it has not been the starting point of constructing a (convergent) numerical method. In a follow-up paper we shall use Langer's approximation functions as ansatz-functions for an $\eps$--uniform numerical method that should also cover turning points.


\section{Conclusion} \label{SEC6}
This paper is concerned with a 1D Schr\"odinger scattering problem in the semi-classical limit, with the inflow given by plane waves. The injection energy and potential are given such that the problem involves both oscillatory and evanescent regions. For the continuous boundary value problem we presented a new, non-overlapping domain decomposition method that separates the oscillatory and evanescent subproblems. The former are treated as IVPs, and the latter as BVPs. Key issues of this approach are the appropriate interface conditions and the final scaling of the solution function. We proved that the domain decomposition method yields the exact solution in a single sweep, performed in the opposite direction of the wave injection.

The hybrid numerical discretization is based on WKB-methods in both types of regions: a WKB-FEM for evanescent regions \cite{Cla}, and a WKB-marching method for oscillatory regions \cite{ABAN11}. The objective of these WKB-methods is to provide an accurate solution -- even on coarse grids and independently of $\epsilon$. Hence, they are \emph{asymptotic preserving}. For the first time we present an error analysis for the WKB-FEM method. Together with the analysis of the WKB-marching method from \cite{ABAN11}, this constitutes the key ingredient for our complete convergence analysis of the hybrid WKB-method. Finally, these error bounds are illustrated and verified in numerical experiments. 

\section{Appendix: proof of Proposition \ref{prop2.2a}}

\
\\
\debproof
\underline{Step 1:} The BVP \eqref{SchBVP} is equivalent to the Schr\"odinger equation \eqref{EQ_ref} on the real line with constant potentials in the leads, and with an incoming plane wave at $x=1$. 
Using the right BC from \eqref{SchBVP}, its solution in the right lead $x\ge1$ hence reads 
\begin{equation}\label{lead-sol}
  \psi(x)=r\,e^{i\frac{\sqrt{a(1)}}{\eps}(x-1)} + e^{-i\frac{\sqrt{a(1)}}{\eps}(x-1)}\,.
\end{equation}
To estimate the reflection coefficient $r$, we consider the current defined in \eqref{current}.
At $x=1$ it reads
$$
  j(1)=\eps  \,\Im [\bar\psi(1)\,\psi'(1)]
  =\eps \sqrt{\alpha(1)} \,(|r|^2-1) \,.
$$
But using the left BC from \eqref{SchBVP} yields
$$
  j(0)=\eps \sqrt{\alpha(0)} \,\Im [\bar \psi(0)\,\psi(0)] = 0\,.
$$
Since the current in a stationary quantum model is constant in $x$, this implies $|r|=1$. Then \eqref{lead-sol} implies 
$$
  |\psi(1)|\le 2\,,\qquad \eps\|\psi'(1)\| \le 2\sqrt{a(1)}\,.
$$
With this bound for the initial condition at $x=1$ we now consider the IVP \eqref{EQ_ref} on $[x_d,1]$. Then, Theorem 2.2 from \cite{Cla} yields the asserted estimate \eqref{2-zone-bound} on the oscillatory region.\\

\noindent
\underline{Step 2:} For the evanescent region $[0,x_d]$ we consider the scaled, real valued solution $\chi$ of the BVP \eqref{region1_ev}. With an elementary argument we first show that $\chi$ has no zero in $[0,x_d]$:\\
Assuming the opposite, let $\chi(x_0)=0$, which implies $\chi'(x_0)\ne0$ (as otherwise $\chi\equiv0$). Then $\chi$ is convex on one side of $x_0$ and concave on the other side, with $\mbox{sgn} (\chi'')=\mbox{sgn} (\chi)$ due to $a\big|_{[0,x_d)}<0$. But then $\chi(0)$ and $\chi'(0)$ would have opposite signs, contradicting the left BC in \eqref{region1_ev}. So we conclude that $\chi$ does not change signs in $[0,x_d]$.\\
Assume now that $\chi(0)<0$ which implies $\chi'(0)<0$ by the left BC in \eqref{region1_ev}. Since then $\chi''<0$ on $[0,x_d)$ we conclude $\chi'(x_d)<0$, contradicting the right BC in \eqref{region1_ev}. \\
This implies that $\chi(0)>0$, and we finally obtain 
$$
  \chi>0\,,\quad \chi'>0\,,\quad \chi''>0\,,\quad \mbox{ on }(0,x_d)\,.
$$
After scaling this auxiliary function, we find that also $|\psi|$ and $|\eps\psi'|$ are increasing on $[0,x_d]$. Therefore the uniform bound \eqref{2-zone-bound} carries over from $x=x_d$ to the evanescent region $[0,x_d]$.
\finproof
\bigskip

\section*{Acknowledgements}
\quad The first author was supported by the FWF-doctoral school ``Dissipation and dispersion in non-linear partial differential equations'', the \"OAD-Amadeus project ``Quantum transport models for semiconductors
and Bose-Einstein condensates'', Universit\'e de Toulouse, and by \emph{Clear Sky Ventures}. He also thanks Kirian D\"opfner for illustrating simulations for \S5. The second author would like to acknowledge support from the ANR project MOONRISE (MOdels, Oscillations and NumeRIcal SchEmes, 2015-2019).

\end{document}